# On Čebotarev Sets

by Kay Wingberg at Heidelberg

The aim of this paper is to define a topology with sufficiently good properties on the set $\mathcal{P}_K$ of prime ideals of a number field $K$. The idea is, roughly speaking, that open sets are given by so-called Čebotarev sets, i.e. sets of the form

$$P_{L|K}(\sigma) = \{\mathfrak{p} \in \mathcal{P}_K \mid \mathfrak{p} \text{ is unramified in } L,\ \sigma = \Big(\frac{L|K}{\mathfrak{P}}\Big),\ \mathfrak{P}|\mathfrak{p}\},$$

where $L|K$ is a finite Galois extension with Galois group $G(L|K)$, $\sigma \in G(L|K)$ and $\big(\frac{L|K}{\mathfrak{P}}\big)$ denotes the Frobenius automorphism with respect to $\mathfrak{P}$, $\mathfrak{P}$ an arbitrary extension of $\mathfrak{p}$ to $L$. The precise definition of the topology $\mathcal{T}_K$ of $\mathcal{P}_K$ is slightly more complicated (see §3) since we want that the natural map

$$\varphi_{K'|K} : (\mathcal{P}_{K'}, \mathcal{T}_{K'}) \longrightarrow (\mathcal{P}_K, \mathcal{T}_K), \quad \mathfrak{P} \mapsto \mathfrak{P} \cap K,$$

is continuous if $K'|K$ is a finite extension. We will show that $(\mathcal{P}_K, \mathcal{T}_K)$ is a strongly zero-dimensional (and so totally disconnected) Hausdorff space with countable base, and so metrizable, hence normal and completely regular (and not discrete). In particular, every point of $(\mathcal{P}_K, \mathcal{T}_K)$ has a base of neighbourhoods consisting of both open and closed sets.

Furthermore we will see that the isolated points of this space have to be prime ideals whose underlying prime numbers ramify in the extension $K|\mathbb{Q}$ (and so the set of isolated points is finite), and that every open neighbourhood of a prime ideal whose underlying prime number is unramified in $K|\mathbb{Q}$ has infinitely many points.

In section 4 we consider uniform structures on $\mathcal{P}_K$ inducing the topology $\mathcal{T}_K$. If $\mathfrak{U}_K$ is the uniformity defined by finite partitions of $\mathcal{P}_K$ given by both open and closed sets, then the completion $(\hat{\mathcal{P}}_K, \hat{\mathfrak{U}}_K)$ of $(\mathcal{P}_K, \mathfrak{U}_K)$ is a profinite space, i.e. compact and totally disconnected.

In section 5 we define an ultra-metric on $\mathcal{P}_K$ inducing the topology $\mathcal{T}_K$ (and a uniform structure which is coarser than $\mathfrak{U}_K$). The idea is that two points $x, y \in \mathcal{P}_K$ are *near*, if they induce in *many* fields with *large* discriminant the same Frobenius automorphism. Again the completion of $\mathcal{P}_K$ with respect to this metric is a profinite space.

# 1 Sets of Completely Decomposed Primes

We start with some remarks on lattices. Let $\mathcal{A}$ and $\mathcal{B}$ be partial ordered sets with respect to the relation $\subseteq$ and let

$$\mathcal{A} \xrightleftharpoons[\psi]{\varphi} \mathcal{B}$$

be maps with the following properties:

I. $\varphi$ and $\psi$ are order-reversing,

II. $A \subseteq \psi\varphi(A)$ and $B \subseteq \varphi\psi(B)$ for all $A \in \mathcal{A}$ and $B \in \mathcal{B}$.

For $A \in \mathcal{A}$ and $B \in \mathcal{B}$ we define

$$\hat{A} := \psi\varphi(A) \text{ and } \hat{B} := \varphi\psi(B),$$

and call $A \in \mathcal{A}$ resp. $B \in \mathcal{B}$ to be saturated if $A = \hat{A}$ resp. $B = \hat{B}$. We set

$$\mathcal{A}_{sat} = \{A \in \mathcal{A} \mid A \text{ is saturated}\}, \qquad \mathcal{B}_{sat} = \{B \in \mathcal{B} \mid B \text{ is saturated}\},$$

and we have the following properties:

(i) $A_1 \subseteq A_2$ implies $\hat{A}_1 \subseteq \hat{A}_2$ and $B_1 \subseteq B_2$ implies $\hat{B}_1 \subseteq \hat{B}_2$.

(ii) $\varphi\psi\varphi = \varphi$, $\psi\varphi\psi = \psi$, $\hat{\hat{A}} = \hat{A}$, $\hat{\hat{B}} = \hat{B}$.

(iii) $\mathcal{B}_{sat}$ is the image of $\mathcal{A}$ under $\varphi$ and $\mathcal{A}_{sat}$ is the image of $\mathcal{B}$ under $\psi$, i.e.

$$\varphi: \mathcal{A} \twoheadrightarrow \mathcal{B}_{sat} \quad \text{and} \quad \psi: \mathcal{B} \twoheadrightarrow \mathcal{A}_{sat}.$$

(iv) $\psi$ and $\varphi$ induce bijections

$$\mathcal{A}_{sat} \xrightleftharpoons[\psi]{\varphi} \mathcal{B}_{sat}.$$

(v) $\quad \varphi(\bigcup_i A_i) = \bigcap_i \varphi(A_i), \quad \psi(\bigcup_i B_i) = \bigcap_i \psi(B_i),$
$\quad \varphi(\bigcap_i \psi(B_i)) = (\bigcup_i B_i)^\wedge, \quad \psi(\bigcap_i \varphi(A_i)) = (\bigcup_i A_i)^\wedge.$

In particular, for saturated sets $A_i$ (resp. $B_i$), $i \in I$, the intersection $\bigcap_i A_i$ (resp. $\bigcap_i B_i$) is saturated.

The verification of these statements is straightforward using the properties I and II of the maps $\varphi$ and $\psi$.

Now let $K$ be a number field and let

$$\mathcal{E}_K = \{L \mid L \text{ is a Galois extension of } K\} \xrightleftharpoons[\psi]{\varphi} \{S \mid S \text{ is a set of primes of } K\} = \mathcal{S}_K$$

where

$\varphi(L) = D(L|K)$ is the set of primes which are completely decomposed in $L|K$,
$\psi(S) = K^S$ is the maximal Galois extension of $K$ which is completely decomposed at $S$.

Obviously, the maps $\varphi$ and $\psi$ are order-reversing and

$$L \subseteq \psi\varphi(L) = K^{D(L|K)} \quad \text{and} \quad S \subseteq \varphi\psi(S) = D(K^S|K).$$

As above we make the following

**Definition 1.1** *For $L \in \mathcal{E}_K$ and $S \in \mathcal{S}_K$ let*

$$\hat{L} = K^{D(L|K)}, \qquad \hat{S} = D(K^S|K).$$

*The extension $L$ is called* **saturated** *if $\hat{L} = L$ and the set $S$ is called* **saturated** *if $\hat{S} = S$.*

From the general remarks above we have bijections

$$(\mathcal{E}_K)_{sat} \underset{\psi}{\overset{\varphi}{\rightleftarrows}} (\mathcal{S}_K)_{sat}$$

and the following

**Lemma 1.2**
(i) $S_1 \subseteq S_2$ implies $\hat{S}_1 \subseteq \hat{S}_2$ and $L_1 \subseteq L_2$ implies $\hat{L}_1 \subseteq \hat{L}_2$.
(ii) $K^S = K^{\hat{S}}$ and $D(L|K) = D(\hat{L}|K)$.
(iii) $D(\prod_i L_i | K) = \bigcap_i D(L_i|K)$ and $K^{\bigcup_i S_i} = \bigcap_i K^{S_i}$.
(iv) $D(\bigcap_i K^{S_i} | K) = (\bigcup_i S_i)^\wedge$ and $K^{\bigcap_i D(L_i|K)} = (\prod_i L_i)^\wedge$.

**Theorem 1.3**
(i) *If $S$ is a finite set of primes, then $S$ is saturated.*
(ii) *If $L$ is a finite Galois extension of $K$, then $L$ is saturated.*

**Proof:** (i) Let $p$ be a prime number and let $L|K$ be a finite Galois extension inside $K^S$. Let $\mathfrak{p}_0 \notin S$, $\mathfrak{P}_0$ a fixed extension of $\mathfrak{p}_0$ to $K^S$ and $\overline{\mathfrak{P}}_0$ the restriction of $\mathfrak{P}_0$ to $L$. By the theorem of Grunwald/Wang (see [4], theorem (9.2.2)) the canonical homomorphism

$$H^1(L, \mathbb{Z}/p\mathbb{Z}) \longrightarrow H^1(L_{\overline{\mathfrak{P}}_0}, \mathbb{Z}/p\mathbb{Z}) \oplus \bigoplus_{\mathfrak{P} \in S(L)} H^1(L_\mathfrak{P}, \mathbb{Z}/p\mathbb{Z})$$

is surjective. In particular, for every $\alpha_{\overline{\mathfrak{P}}_0} \in H^1(L_{\overline{\mathfrak{P}}_0}, \mathbb{Z}/p\mathbb{Z}))$ there exists an element $\beta \in H^1(L, \mathbb{Z}/p\mathbb{Z})$ which is mapped to $(\alpha_{\overline{\mathfrak{P}}_0}, 0, \ldots, 0)$. But $\beta$ lies in the subgroup $H^1(K^S|L, \mathbb{Z}/p\mathbb{Z})$ of $H^1(L, \mathbb{Z}/p\mathbb{Z})$. Therefore

$$H^1(K^S|L, \mathbb{Z}/p\mathbb{Z}) \longrightarrow H^1(L_{\overline{\mathfrak{P}}_0}, \mathbb{Z}/p\mathbb{Z})$$

is surjective. It follows that the completion of $K^S$ with respect to the prime $\mathfrak{P}_0$, $\mathfrak{P}_0|\mathfrak{p}_0$ and $\mathfrak{p}_0 \notin S$, is equal to the algebraic closure of $K_{\mathfrak{p}_0}$ (since $G(\overline{K}_{\mathfrak{p}_0}|K_{\mathfrak{p}_0})$ is pro-solvable). In particular, $D(K^S|K) = S$.

(ii) Let $L'$ be a finite Galois extension of $K$ with $L \subseteq L' \subseteq K^{D(L|K)}$. Since $D(L|K) \subseteq D(L'|K)$, we obtain for the densities of these sets the inequality $\delta(D(L|K)) \leq \delta(D(L'|K))$, and so, by Čebotarev's density theorem,

$$[L' : K] = \delta(D(L'|K))^{-1} \leq \delta(D(L|K))^{-1} = [L : K].$$

This shows that $L' = L$ and so $L = K^{D(L|K)}$. $\square$

## 2  Čebotarev Sets

Let $K$ be a number field and let $\mathcal{P}_K$ be the set of all prime ideals $\mathfrak{p} \neq (0)$ of $K$. For a finite Galois extension $L|K$ with Galois group $G(L|K)$ we denote the set of prime ideals of $K$ which are unramified in $L$ by $U(L|K)$ and the set of ramified prime ideals by $R(L|K)$. For an element $\sigma \in G(L|K)$ let

$$P_{L|K}(\sigma) = \{\mathfrak{p} \in U(L|K) \,|\, \sigma = \Big(\frac{L|K}{\mathfrak{P}}\Big) \text{ for a prime ideal } \mathfrak{P}|\mathfrak{p} \text{ of } L\},$$

where $\Big(\frac{L|K}{\mathfrak{P}}\Big)$ denotes the Frobenius automorphism with respect to $\mathfrak{P}$. Obviously, this set depends only on the conjugacy class $\langle\!\langle\sigma\rangle\!\rangle = \{\tau\sigma\tau^{-1} \,|\, \tau \in G(L|K)\}$ of $\sigma$ and

$$P_{L|K}(\sigma) \cap P_{L|K}(\tau) = \varnothing \quad \text{if } \langle\!\langle\sigma\rangle\!\rangle \neq \langle\!\langle\tau\rangle\!\rangle,$$

and $P_{L|K}(1) = D(L|K)$. By Čebotarev's density theorem we have

$$\delta(P_{L|K}(\sigma)) = \frac{\#\langle\!\langle\sigma\rangle\!\rangle}{\#G(L|K)},$$

where $\delta(S) = \delta_K(S)$ denotes the Dirichlet density of a set $S$ of primes of $K$. Observe that for a finite Galois extension $L|K$ and a set $S(K)$ of primes of $K$ we have

$$\delta_L(S(L)) = \delta_K(S(K) \cap D(L|K))) \cdot [L : K],$$

where $S(L)$ denotes the set of all extensions of $S(K)$ to $L$. For sets $S_1$ and $S_2$ of primes, we use the notation

$$S_1 \subsetsim S_2 :\Longleftrightarrow \delta(S_1 \setminus S_2) = 0,$$

i.e. $S_1$ is contained in $S_2$ up to a set of primes of density zero, and

$$S_1 \stackrel{.}{=} S_2 :\Longleftrightarrow S_1 \subsetsim S_2 \text{ and } S_2 \subsetsim S_1.$$

**Proposition 2.1** Let $N|L|K$ be finite Galois extensions, $H = G(N|L)$ and let $\bar{\sigma} \in G(L|K)$. Then

$$U(N|K) \cap P_{L|K}(\bar{\sigma}) = \bigcup_{\langle\!\langle\tau\rangle\!\rangle \cap \sigma H \neq \emptyset}^{\cdot} P_{N|K}(\tau),$$

where $\sigma$ is a lifting of $\bar{\sigma}$ to $G(N|K)$; in particular

$$U(N|K) \cap D(L|K) = \bigcup_{\langle\!\langle\tau\rangle\!\rangle \cap H \neq \emptyset}^{\cdot} P_{N|K}(\tau).$$

**Proof:** Let $\mathfrak{p}$ be a prime ideal of $K$ which is unramified in $N|K$. Then $\mathfrak{p} \in P_{L|K}(\bar{\sigma})$ if and only if there exists a prime $\overline{\mathfrak{P}}|\mathfrak{p}$ of $L$ such that $\bar{\sigma} = \left(\frac{L|K}{\overline{\mathfrak{P}}}\right)$, i.e. if there exists a prime $\mathfrak{P}|\mathfrak{p}$ of $N$ such that $\sigma H = \left(\frac{N|K}{\mathfrak{P}}\right) H$. This ist equivalent to the assertion that there exists an element in $\sigma H$ which is contained in the conjugacy class $\langle\!\langle\tau\rangle\!\rangle$ of $\tau = \left(\frac{N|K}{\mathfrak{P}}\right)$ for some prime ideal $\mathfrak{P}|\mathfrak{p}$ of $N$, i.e. if $\mathfrak{p} \in P_{N|K}(\tau)$ for some $\tau \in G(N|K)$ with $\langle\!\langle\tau\rangle\!\rangle \cap \sigma H \neq \emptyset$. $\square$

**Corollary 2.2** Let $L_1|K$ and $L_2|K$ be finite Galois extensions, $H_i = G(L_1L_2|L_i)$ and $\bar{\sigma}_i \in G(L_i|K)$, $i = 1, 2$. Then

$$P_{L_1|K}(\bar{\sigma}_1) \cap P_{L_2|K}(\bar{\sigma}_2) = \bigcup_{\substack{\langle\!\langle\tau\rangle\!\rangle \cap \sigma_1 H_1 \neq \emptyset \\ \langle\!\langle\tau\rangle\!\rangle \cap \sigma_2 H_2 \neq \emptyset}}^{\cdot} P_{L_1L_2|K}(\tau)$$

**Proof:** This is a consequence of proposition 2.1 since $U(L_1|K) \cap U(L_2|K) = U(L_1L_2|K)$. $\square$

If $\bar{\sigma}_1 = 1 = \bar{\sigma}_2$, then the corollary above is just the assertion $D(L_1|K) \cap D(L_2|K) = D(L_1L_2|K)$ used in section 1.

**Corollary 2.3** Let $L_1|K$ and $L_2|K$ be finite Galois extensions and $\sigma_i \in G(L_i|K)$, $i = 1, 2$. Then

$$P_{L_1|K}(\sigma_1) \cap P_{L_2|K}(\sigma_2) \neq \emptyset$$

if and only if

$$(\sigma_1)^{-1}(\sigma_2)^\rho \in G(L_1L_2|L_1 \cap L_2) \quad \text{for some } \rho \in G(L_1L_2|K)$$

(here $\sigma_i$ denotes an arbitrary lifting of $\sigma_i$ to $G(L_1L_2|K)$). In particular, if $L_1$ and $L_2$ are linearly disjoint over $K$, then all sets $P_{L_1|K}(\sigma_1)$ and $P_{L_2|K}(\sigma_2)$ have a non-trivial intersection.

For an element $\tau$ of a finite group $G$ we denote the stabilizer of $\tau$ under conjugation by $G$ by $\mathrm{St}_G(\tau)$.

**Proposition 2.4** *Let $L_1$ and $L_2$ be finite Galois extensions of $K$. For an element $\sigma_i \in G(L_i|K)$ we denote its reduction modulo $G(L_i|L_1 \cap L_2)$ by $\overline{\sigma_i}$, $i = 1, 2$. Then the following assertions are equivalent:*

(i) $P_{L_1|K}(\sigma_1) \subsetneq P_{L_2|K}(\sigma_2)$,

(ii) $\langle\!\langle \overline{\sigma_1} \rangle\!\rangle = \langle\!\langle \overline{\sigma_2} \rangle\!\rangle$ *and* $\#\mathrm{St}_{G(L_2|K)}(\sigma_2) = \#\mathrm{St}_{G(L_1 \cap L_2|K)}(\overline{\sigma_2})$.

*In particular, $P_{L_1|K}(\sigma_1) \doteq P_{L_2|K}(\sigma_2)$ if and only if $\langle\!\langle \overline{\sigma_1} \rangle\!\rangle = \langle\!\langle \overline{\sigma_2} \rangle\!\rangle$ and*

$$\#\mathrm{St}_{G(L_1|K)}(\sigma_1) = \#\mathrm{St}_{G(L_1\cap L_2|K)}(\overline{\sigma_1}) = \#\mathrm{St}_{G(L_2|K)}(\sigma_2).$$

**Proof:** Let $N = L_1 L_2$, $H_i = G(N|L_i)$, $i = 1, 2$, and $H = G(L_2|L_1 \cap L_2) \cong H_1$. We lift $\sigma_i$ to $G(N|K)$ and denote it again by $\sigma_i$. Assume that (i) holds, i.e.

$$P_{L_1|K}(\sigma_1) \doteq \bigcup_{\langle\!\langle \tau \rangle\!\rangle \cap \sigma_1 H_1 \neq \varnothing} P_{N|K}(\tau) \subsetneq \bigcup_{\langle\!\langle \tau \rangle\!\rangle \cap \sigma_2 H_2 \neq \varnothing} P_{N|K}(\tau) \doteq P_{L_2|K}(\sigma_2).$$

Since the sets $P_{N|K}(\tau)$ have positive density, it follows that for every $h_1 \in H_1$ there exist $h_2 \in H_2$ and $\rho \in G(N|K)$ such that $\sigma_1 h_1 = (\sigma_2)^\rho h_2$, and so $\langle\!\langle \overline{\sigma_1} \rangle\!\rangle = \langle\!\langle \overline{\sigma_2} \rangle\!\rangle$.

If $h \in H$ is a fixed element and $\tilde{h}$ a lifting of $h$ to $G(N|K)$, then it follows that for every $h_1 \in H_1$ there exist $h_2 \in H_2$ and $\rho \in G(N|K)$ such that $\sigma_1 h_1 = (\sigma_2 \tilde{h})^\rho h_2$, and so

$$\bigcup_{\langle\!\langle \tau \rangle\!\rangle \cap \sigma_1 H_1 \neq \varnothing} P_{N|K}(\tau) \subsetneq \bigcup_{\langle\!\langle \tau \rangle\!\rangle \cap (\sigma_2 \tilde{h}) H_2 \neq \varnothing} P_{N|K}(\tau).$$

We obtain $P_{L_1|K}(\sigma_1) \subsetneq P_{L_2|K}(\sigma_2 h)$ for every $h \in H$, i.e.

$$P_{L_1|K}(\sigma_1) \subsetneq \bigcap_{h \in H} P_{L_2|K}(\sigma_2 h).$$

Since

$$P_{L_2|K}(\sigma_2) \cap P_{L_2|K}(\sigma_2 h) \neq \varnothing \quad \text{if and only if} \quad \langle\!\langle \sigma_2 \rangle\!\rangle = \langle\!\langle \sigma_2 h \rangle\!\rangle,$$

it follows that

$$P_{L_1 \cap L_2|K}(\sigma_2 H) \doteq \bigcup_{\langle\!\langle \sigma_2 h \rangle\!\rangle, h \in H} P_{L_2|K}(\sigma_2 h) = P_{L_2|K}(\sigma_2),$$

and therefore

$$\frac{\#\langle\!\langle \sigma_2 H \rangle\!\rangle_{G(L_1 \cap L_2|K)}}{\#G(L_1 \cap L_2|K)} = \frac{\#\langle\!\langle \sigma_2 \rangle\!\rangle_{G(L_2|K)}}{\#G(L_2|K)}.$$

From this equation we get

$$\#\mathrm{St}_{G(L_2|K)}(\sigma_2) = \#\mathrm{St}_{G(L_1 \cap L_2|K)}(\sigma_2 H).$$

Conversely, using the arguments above in the other direction, we obtain from the assertion (ii) that

$$P_{L_1 \cap L_2|K}(\sigma_2 H) \eqsim P_{L_2|K}(\sigma_2).$$

Since $\langle\langle \overline{\sigma_1} \rangle\rangle = \langle\langle \overline{\sigma_2} \rangle\rangle$, we get

$$P_{L_1|K}(\sigma_1) \subsetneq P_{L_1 \cap L_2|K}(\sigma_2 H)$$

and so (i). This finishes the proof of the proposition. $\square$

**Remark:** From the proposition above it follows that it is possible that

$$P_{L|K}(\overline{\sigma}) \eqsim P_{N|K}(\sigma),$$

where $L$ is a *proper* subfield of the extension $N|K$, $\sigma$ is an element of $G(N|K)$ and $\overline{\sigma} = \sigma G(N|L)$. As an example let $G(N|K)$ be the non-abelian group of order $p^3$, $p$ an odd prime number, with generators $\sigma, \tau, \rho$ and the defining relations

$$\sigma^p = \tau^p = \rho^p = [\sigma, \rho] = [\tau, \rho] = 1 \quad \text{and} \quad [\sigma, \tau] = \rho,$$

and let $G(N|L)$ be the normal subgroup generated by $\rho$.

The following corollary is a generalization of a theorem of M.Bauer (see [3], theorem (13.9)).

**Corollary 2.5** *Let $L_1$ and $L_2$ be finite Galois extensions of $K$ and let $\sigma_i$ be an element of $G(L_i|K)$, $i = 1, 2$. Assume that $\sigma_2$ lies in the center of $G(L_2|K)$. Then the following assertions are equivalent:*

(i) $P_{L_1|K}(\sigma_1) \subsetneq P_{L_2|K}(\sigma_2)$,
(ii) $L_2 \subseteq L_1$ and $\sigma_2 = \sigma_1 \bmod G(L_1|L_2)$.

*In particular, if $\sigma_i$ lies in the center of $G(L_i|K)$, $i = 1, 2$, then*

$$P_{L_1|K}(\sigma_1) \eqsim P_{L_2|K}(\sigma_2) \quad \textit{if and only if} \quad L_1 = L_2 \textit{ and } \sigma_1 = \sigma_2.$$

**Proof:** By assumption $\sigma_2$ lies in the center of $G(L_2|K)$, and so $\#\mathrm{St}_{G(L_2|K)}(\sigma_2) = \#\mathrm{St}_{G(L_1 \cap L_2|K)}(\overline{\sigma_2})$ if and only if $G(L_2|L_1 \cap L_2) = 1$, i.e. $L_2 \subseteq L_1$. Now the corollary follows from proposition 2.4. $\square$

**Application to quadratic fields**

We will use the following notation: the non-trivial element of the Galois group $G(\mathbb{Q}(\sqrt{a})|\mathbb{Q})$, $a$ a squarefree integer, is denoted by $-1$. If $\sigma \in \{\pm 1\}$, then $\sigma' \in \{\pm 1\}$ is defined by $\sigma\sigma' = -1$. For $x, y \in \mathcal{P}_\mathbb{Q}$, $y$ odd, $\left(\frac{x}{y}\right)$ is the Legendre symbol. Let $\varepsilon(\ell) = (\ell - 1)/2$ if $\ell$ is odd, and zero otherwise, and let

$$L_p = \mathbb{Q}(\sqrt{(-1)^{\varepsilon(p)}p}), \quad p \in \mathcal{P}_\mathbb{Q}.$$

For $\tilde{\sigma} = (\sigma_p)_p \in \prod_{p \in \mathcal{P}_\mathbb{Q}} \{\pm 1\}$ and $\varepsilon \in \{\pm 1\}$ let

$$\delta^\varepsilon(\tilde{\sigma}) = \delta(\{p \in \mathcal{P}_\mathbb{Q} \mid \sigma_p = \varepsilon\})$$

be the density of the set of prime numbers $p$ such that $\sigma_p$ has a given value.

**Proposition 2.6** *Let $\tilde{\sigma} = (\sigma_p)_p \in \prod_{p \in \mathcal{P}_\mathbb{Q}} \{\pm 1\}$ and $S \subseteq \mathcal{P}_\mathbb{Q}$ a set of density equal to 1. Then*

(i) $$\#\left(\mathcal{P}_\mathbb{Q} \setminus \bigcup_{p \in S} P_{L_p | \mathbb{Q}}(\sigma_p)\right) \leq 1,$$

(ii) *if* $\#\left(\mathcal{P}_\mathbb{Q} \setminus \bigcup_{p \in S} P_{L_p | \mathbb{Q}}(\sigma_p)\right) = 1$, *then* $\delta^1(\tilde{\sigma}) = \delta^{-1}(\tilde{\sigma}) = 1/2$.

**Proof:** Let

$$q_1, q_2 \in \bigcap_{p \in S} \left(P_{L_p | \mathbb{Q}}(\sigma'_p) \cup \{p\}\right).$$

If $q_1$ and $q_2$ are odd, then

$$\left(\frac{(-1)^{\varepsilon(p)}p}{q_1}\right) = \sigma'_p = \left(\frac{(-1)^{\varepsilon(p)}p}{q_2}\right), \quad \text{and so} \quad \left(\frac{q_1}{p}\right) = \left(\frac{q_2}{p}\right)$$

for all $p \in S \setminus \{2, q_1, q_2\}$. From corollary 2.5 it follows that $q_1 = q_2$. If $q_2 = 2$, then

$$2 \in P_{L_p | \mathbb{Q}}(\sigma'_p), \ p \text{ odd}, \quad \text{if and only if} \quad \left(\frac{2}{p}\right) = \sigma'_p$$

($(-1)^{\varepsilon(p)}p \equiv 1 \mod 4$ and 2 splits in $L_p$ if and only if $(-1)^{\varepsilon(p)}p \equiv 1 \mod 8$, i.e. $\left(\frac{2}{p}\right) = 1$.) Again it follows that $\left(\frac{q_1}{p}\right) = \left(\frac{2}{p}\right)$ for all $p \in S \setminus \{2, q_1\}$, and so $q_1 = 2$. Since

$$\mathcal{P}_\mathbb{Q} \setminus \bigcup_{p \in S} P_{L_p | \mathbb{Q}}(\sigma_p) = \bigcap_{p \in S} \left(P_{L_p | \mathbb{Q}}(\sigma'_p) \cup \{p\}\right),$$

we proved (i).

If $q \in \bigcap_{p \in S} \left(P_{L_p | \mathbb{Q}}(\sigma'_p) \cup \{p\}\right)$, and so $q \in \bigcap_{p \in S, p \neq q} P_{L_p | \mathbb{Q}}(\sigma'_p)$, then $\left(\frac{q}{p}\right) = \sigma'_p$ for all $p \in S \setminus \{2, q\}$. It follows that $\delta^1(\tilde{\sigma}) = \delta^{-1}(\tilde{\sigma}) = 1/2$. □

**Definition 2.7** *A set $S$ of prime ideals of $K$ is called **Čebotarev set** if there exist a finite Galois extension $L$ of $K$ and an element $\sigma \in G(L|K)$ such that*
$$S = P_{L|K}(\sigma).$$
*We set*
$$\mathcal{C}_K = \{S \subseteq \mathcal{P}_K \text{ is a Čebotarev set}\}.$$

For a finite extension $K'|K$ let
$$\varphi_{K'|K} : \mathcal{P}_{K'} \longrightarrow \mathcal{P}_K, \quad \mathfrak{P} \mapsto \mathfrak{p} = \mathfrak{P} \cap K,$$
and we also denote the corresponding map on the set of all subsets of $\mathcal{P}_{K'}$ by $\varphi_{K'|K}$. For a subfield $E \subseteq K$, a finite Galois extension $F|E$ and an element $\sigma \in G(F|E)$,

$$\begin{array}{ccc} K & & F \\ | & \diagup & \\ & \text{\scriptsize Galois} & \\ E & & \end{array}$$

let
$$P_{F|E}^K(\sigma) = \varphi_{K|E}^{-1} P_{F|E}(\sigma),$$
$$U^K(F|E) = \varphi_{K|E}^{-1} U(F|E) \quad \text{and} \quad R^K(F|E) = \varphi_{K|E}^{-1} R(F|E).$$

**Definition 2.8** *For a number field $K$ let*
$$\mathcal{B}_K = \bigcup_{E \subseteq K} \varphi_{K|E}^{-1} \mathcal{C}_E = \{P_{F|E}^K(\sigma) \,|\, E \subseteq K, \, F|E \text{ finite Galois}, \, \sigma \in G(F|E)\}$$
*and*
$$\mathcal{A}_K = \{\bigcap_{i=1}^n P_{F_i|E_i}^K(\sigma_i) \,|\, n \in \mathbb{N}, \, P_{F_i|E_i}^K(\sigma_i) \in \mathcal{B}_K\}.$$
*Observe that $\mathcal{C}_K \subseteq \mathcal{B}_K \subseteq \mathcal{A}_K$.*

In the next section we will need the following two lemmas.

**Lemma 2.9** *Let $K|\mathbb{Q}$ be a finite Galois extension and for $i = 1, \ldots, n$ let $E_i \subseteq K$ be subfields of $K$, $F_i|E_i$ finite Galois extensions and $\sigma_i \in G(F_i|E_i)$. Then*
$$P_{K|\mathbb{Q}}^K(1) \cap \bigcap_{i=1}^n P_{F_i|E_i}^K(\sigma_i)$$
*is empty or has positive density.*

**Proof:** Let $F|E$ be one of the extensions $F_i|E_i$. Then

$$\varphi_{K|\mathbb{Q}}^{-1}P_{K|\mathbb{Q}}(1) \cap \varphi_{K|E}^{-1}P_{F|E}(\sigma) = \varphi_{K|\mathbb{Q}}^{-1}P_{K|\mathbb{Q}}(1) \cap \bigcup_{\langle\!\langle\tau\rangle\!\rangle \cap \sigma H \neq \emptyset} P_{FK|K}(\tau),$$

where $H = G(FK|F)$. Indeed, let $\mathfrak{P}_K \in \varphi_{K|\mathbb{Q}}^{-1}P_{K|\mathbb{Q}}(1) \subseteq \mathcal{P}_K$ and consider the diagram of fields

```
        FK
       /  \
      F    K
       \   |
       F∩K
         |
         E
```

Let $\mathfrak{P}_F$ be an extension of $\mathfrak{p} = \mathfrak{P}_K \cap E$ to $F$ and $\mathfrak{P}_{FK}$ an extension of $\mathfrak{P}_F$ to $FK$. Then $\mathfrak{P}'_K = \mathfrak{P}_{FK} \cap K$ is conjugated to $\mathfrak{P}_K$. Since $\mathfrak{P}_{FK}$ is unramified over $K$ and the residue degree $f(\mathfrak{P}'_K|\mathfrak{p}) = 1$, we have

$$\left(\frac{FK|K}{\mathfrak{P}_{FK}}\right)_{|F} = \left(\frac{F|E}{\mathfrak{P}_F}\right).$$

Now the equality stated above follows easily. Thus we obtain

$$P_{K|\mathbb{Q}}^K(1) \cap \bigcap_{i=1}^n P_{F_i|E_i}^K(\sigma_i)$$
$$= \varphi_{K|\mathbb{Q}}^{-1}P_{K|\mathbb{Q}}(1) \cap \bigcap_{i=1}^n \bigcup_{\langle\!\langle\tau_i\rangle\!\rangle \cap \sigma_i H_i \neq \emptyset} P_{F_iK|K}(\tau_i)$$
$$= \varphi_{K|\mathbb{Q}}^{-1}P_{K|\mathbb{Q}}(1) \cap \bigcup_{\langle\!\langle\tau_1\rangle\!\rangle \cap \sigma_1 H_1 \neq \emptyset} \cdots \bigcup_{\langle\!\langle\tau_n\rangle\!\rangle \cap \sigma_n H_n \neq \emptyset} \Big(P_{F_1K|K}(\tau_1) \cap \cdots \cap P_{F_nK|K}(\tau_n)\Big).$$

From corollary 2.2 it follows that the sets

$$P_{F_1K|K}(\tau_1) \cap \cdots \cap P_{F_nK|K}(\tau_n)$$

are empty or have positive density. Since the density of $\varphi_{K|\mathbb{Q}}^{-1}P_{K|\mathbb{Q}}(1)$ is equal to 1, we proved the lemma. $\square$

**Lemma 2.10** *Let $K$ be a number field and for $i = 1, \ldots, n$ let $E_i \subseteq K$ be subfields of $K$, $F_i|E_i$ finite Galois extensions, $E = \bigcap_i E_i$ and $\sigma_i \in G(F_i|E_i)$. Then the set*

$$S = \varphi_{K|E}^{-1}U(K|E) \cap \bigcap_{i=1}^n P_{F_i|E_i}^K(\sigma_i)$$

*is empty or infinite.*

**Proof:** We may assume that $K = F_1 = \cdots = F_n$ and that $K|E$ is a Galois extension. In order to see this let $M$ be the normal closure of $KF_1 \cdots F_n$ over $E$ and $S_M = \varphi_{M|K}^{-1}(S)$. Then $S_M$ is non-empty and finite if and only $S$ is. We get

$$\begin{aligned}
S_M &= \varphi_{M|E}^{-1} U(K|E) \cap \varphi_{M|E_1}^{-1} P_{F_1|E_1}(\sigma_1) \cap \cdots \cap \varphi_{M|E_n}^{-1} P_{F_n|E_n}(\sigma_n) \\
&= \varphi_{M|E}^{-1} U(M|E) \cap \varphi_{M|E_1}^{-1} P_{F_1|E_1}(\sigma_1) \cap \cdots \cap \varphi_{M|E_n}^{-1} P_{F_n|E_n}(\sigma_n) \\
&= \varphi_{M|E}^{-1} U(M|E) \cap \bigcup_{\langle\!\langle \tau_1 \rangle\!\rangle \cap \sigma_1 H_1 \neq \varnothing} \varphi_{M|E_1}^{-1} P_{M|E_1}(\tau_1) \cap \cdots \cap \bigcup_{\langle\!\langle \tau_n \rangle\!\rangle \cap \sigma_n H_n \neq \varnothing} \varphi_{M|E_n}^{-1} P_{M|E_n}(\tau_n) \\
&= \varphi_{M|E}^{-1} U(M|E) \cap \bigcup_{\substack{\langle\!\langle \tau_i \rangle\!\rangle \cap \sigma_i H_i \neq \varnothing \\ i=1,\ldots,n}} \left( \varphi_{M|E_1}^{-1} P_{M|E_1}(\tau_1) \cap \cdots \cap \varphi_{M|E_n}^{-1} P_{M|E_n}(\tau_n) \right),
\end{aligned}$$

where $H_i = G(M|F_i)$, $i = 1, \ldots, n$, and so it is enough to show that the sets

$$S(\tau_1, \ldots, \tau_n) = \varphi_{M|E}^{-1} U(M|E) \cap \varphi_{M|E_1}^{-1} P_{M|E_1}(\tau_1) \cap \cdots \cap \varphi_{M|E_n}^{-1} P_{M|E_n}(\tau_n)$$

are empty or infinite. For a set $T$ of primes of $M$ let $(T)_{G(M|E)}$ be the closure under conjugation by $G(M|E)$. Obviously, it is sufficient to show that $(S(\tau_1, \ldots, \tau_n))_{G(M|E)}$ is empty or infinite.

Suppose that $\mathfrak{P} \in S(\tau_1, \ldots, \tau_n)$ and let $\mathfrak{p} = \mathfrak{P} \cap E$. We fix an index $i$. Then $\mathfrak{P}_{E_i} = \mathfrak{P} \cap E_i \in P_{M|E_i}(\tau_i)$, i.e. there exists an extension $\mathfrak{P}_M$ of $\mathfrak{P}_{E_i}$ in $M$ such that $\tau_i = \left( \frac{M|E_i}{\mathfrak{P}_M} \right)$. Since $\mathfrak{P}$ and $\mathfrak{P}_M$ are conjugated over $E_i$, it follows that there is an element $\rho_i \in G(M|E_i)$ such that $\tau_i^{\rho_i} = \left( \frac{M|E_i}{\mathfrak{P}} \right)$ and we may assume that $\tau_i = \left( \frac{M|E_i}{\mathfrak{P}} \right)$.

Since $\mathfrak{P} \in \varphi_{M|E}^{-1} U(M|E)$, we have the element $\sigma = \left( \frac{M|E}{\mathfrak{P}} \right) \in G(M|E)$ and it follows that

$$\tau_i = \left( \frac{M|E_i}{\mathfrak{P}} \right) = \left( \frac{M|E}{\mathfrak{P}} \right)^{f(\mathfrak{P}_{E_i}|\mathfrak{p})} = \sigma^{f(\mathfrak{P}_{E_i}|\mathfrak{p})},$$

where $f(\mathfrak{P}_{E_i}|\mathfrak{p})$ is the inertia degree of $\mathfrak{P}_{E_i}$ over $E$. We claim that

$$\varphi_{M|E}^{-1} P_{M|E}(\sigma) \subseteq (S(\tau_1, \ldots, \tau_n))_{G(M|E)}.$$

Indeed, let $\mathfrak{P}' \in \varphi_{M|E}^{-1} P_{M|E}\left( \left( \frac{M|E}{\mathfrak{P}} \right) \right)$. Then there exists a prime $\mathfrak{P}''$ of $M$ which is conjugated to $\mathfrak{P}'$ over $E$ such that $\left( \frac{M|E}{\mathfrak{P}} \right) = \left( \frac{M|E}{\mathfrak{P}''} \right)$. Let $f(\mathfrak{P}''_{E_i}|\mathfrak{p})$ be the inertia degree of $\mathfrak{P}''_{E_i}$ over $E$. Since

$$G(M|E_i) \cap G_{\mathfrak{P}''}(M|E) = G_{\mathfrak{P}''}(M|E_i),$$

we get

$$\Big(\frac{M|E}{\mathfrak{P}''}\Big)^{f(\mathfrak{P}_{E_i}|\mathfrak{p})} = \Big(\frac{M|E}{\mathfrak{P}}\Big)^{f(\mathfrak{P}_{E_i}|\mathfrak{p})} = \Big(\frac{M|E_i}{\mathfrak{P}}\Big) \in G_{\mathfrak{P}''}(M|E_i).$$

Since $G_{\mathfrak{P}''}(M|E_i)$ is generated by $\Big(\frac{M|E_i}{\mathfrak{P}''}\Big) = \Big(\frac{M|E}{\mathfrak{P}''}\Big)^{f(\mathfrak{P}''_{E_i}|\mathfrak{p})}$, $f(\mathfrak{P}''_{E_i}|\mathfrak{p})$ divides $f(\mathfrak{P}_{E_i}|\mathfrak{p})$. Analogously,

$$\Big(\frac{M|E}{\mathfrak{P}}\Big)^{f(\mathfrak{P}''_{E_i}|\mathfrak{p})} = \Big(\frac{M|E}{\mathfrak{P}''}\Big)^{f(\mathfrak{P}''_{E_i}|\mathfrak{p})} = \Big(\frac{M|E_i}{\mathfrak{P}''}\Big) \in G_{\mathfrak{P}}(M|E_i),$$

and so $f(\mathfrak{P}_{E_i}|\mathfrak{p})$ divides $f(\mathfrak{P}''_{E_i}|\mathfrak{p})$. Therefore we obtain

$$\Big(\frac{M|E_i}{\mathfrak{P}''}\Big) = \Big(\frac{M|E}{\mathfrak{P}''}\Big)^{f(\mathfrak{P}''_{E_i}|\mathfrak{p})} = \Big(\frac{M|E}{\mathfrak{P}}\Big)^{f(\mathfrak{P}_{E_i}|\mathfrak{p})} = \Big(\frac{M|E_i}{\mathfrak{P}}\Big) = \tau_i.$$

It follows that $\mathfrak{P}'' \in \varphi_{M|E_i}^{-1} P_{M|E_i}(\tau_i)$ for all $i = 1, \ldots, n$, i.e. $\mathfrak{P}'' \in S(\tau_1, \ldots, \tau_n)$, and so $\mathfrak{P}' \in (S(\tau_1, \ldots, \tau_n))_{G(M|E)}$. This proves the claim. Since $\varphi_{M|E}^{-1} P_{M|E}(\sigma)$ is an infinite set, we proved the lemma. □

We finish this section with a slightly more general version of the theorem of Grunwald/Wang (see also [4], theorem (9.2.2)).

Let $p$ be a prime number, $K$ a number field and $S \supseteq T$ sets of primes of $K$, where $S$ contains the set $S_p \cup S_\infty$ of archimedean primes and primes above $p$. Let $K_S$ be the maximal extension of $K$ which is unramified outside $S$. By $\mu_p$ we denote the group of all $p$-th roots of unity.

**Theorem 2.11** *Let $K$ be a number field and let $S \supseteq T$ be sets of primes of $K$, where $S \supseteq S_p \cup S_\infty$, $T$ is finite and*

$$\delta(S \cap D(K(\mu_p)|K)) > \frac{1}{p\,[K(\mu_p):K]}.$$

*Then the canonical homomorphism*

$$H^1(K_S|K, \mathbb{Z}/p\mathbb{Z}) \longrightarrow \bigoplus_{\mathfrak{p} \in T} H^1(K_\mathfrak{p}, \mathbb{Z}/p\mathbb{Z})$$

*is surjective.*

**Proof:** Using [4], lemma (9.2.1), it is enough to show that the canonical map

$$H^1(K_S|K, \mu_p) \longrightarrow \prod_{\mathfrak{p} \in (S \setminus T)(K)} H^1(K_\mathfrak{p}, \mu_p)$$

is injective. Since $[K(\mu_p) : K]$ is prime to $p$, it is sufficient to show the injectivity of the homomorphism

$$H^1(K_S|K(\mu_p), \mu_p) \longrightarrow \prod_{\mathfrak{P} \in (S \setminus T)(K(\mu_p))} H^1(K(\mu_p)_{\mathfrak{P}}, \mu_p).$$

An element of the kernel corresponds to a Galois extension $L|K(\mu_p)$ of degree $p$ which is unramified outside $S(K(\mu_p))$ and completely decomposed at $(S \setminus T)(K(\mu_p))$. Since

$$\begin{aligned}\delta_{K(\mu_p)}((S \setminus T)(K(\mu_p))) &= \delta_{K(\mu_p)}(S(K(\mu_p)) \\ &= \delta_K(S(K) \cap D(K(\mu_p)|K)) \cdot [K(\mu_p) : K] > \tfrac{1}{p},\end{aligned}$$

such an extension has to be trivial. $\square$

# 3 Topology

In this section we define a topology on the set $\mathcal{P}_K$ of non-trivial prime ideals of a number field $K$.

**Definition 3.1** *For a number field $K$ let $\mathcal{T}_K$ be the topology on $\mathcal{P}_K$ generated by*

$$\mathcal{A}_K = \{\bigcap_{i=1}^{n} P^K_{F_i|E_i}(\sigma_i) \mid n \in \mathbb{N},\ P^K_{F_i|E_i}(\sigma_i) \in \mathcal{B}_K\},$$

*i.e. $\mathcal{A}_K$ is a base and*

$$\mathcal{B}_K = \{P^K_{F|E}(\sigma) \mid E \subseteq K,\ F|E \text{ a finite Galois extension},\ \sigma \in G(F|E)\}$$

*is a subbase of $\mathcal{T}_K$. Obviously, the topology $\mathcal{T}_K$ has a countable base.*

**Remark:** From proposition 2.2 it follows that

$$\mathcal{C}_{\mathbb{Q}} \cup \{\varnothing\} = \{P_{F|\mathbb{Q}}(\sigma) \mid F|\mathbb{Q} \text{ a finite Galois extension},\ \sigma \in G(F|\mathbb{Q})\} \cup \{\varnothing\}$$

is a base of $\mathcal{T}_{\mathbb{Q}}$.

**Proposition 3.2** *If $K'|K$ is a finite extension, then the map*

$$\varphi_{K'|K} : (\mathcal{P}_{K'}, \mathcal{T}_{K'}) \longrightarrow (\mathcal{P}_K, \mathcal{T}_K), \quad \mathfrak{P} \mapsto \mathfrak{P} \cap K,$$

*is continuous.*

**Proof:** This follows directly from the definition of the topologies $\mathcal{T}_K$ and $\mathcal{T}_{K'}$. $\square$

An observation, which is not quite obvious, is the following:

$$\mathcal{T}_K \text{ is not the discrete topology on } \mathcal{P}_K.$$

In order to see this, let us first assume that $K|\mathbb{Q}$ is a Galois extension and suppose that $\mathcal{T}_K$ is the discrete topology. Then for every point $\mathfrak{p} \in \mathcal{P}_K$ the set $\{\mathfrak{p}\}$ is open and therefore there exist finite Galois extensions $F_i|E_i$, $E_i \subseteq K$, and $\sigma_i \in G(F_i|E_i)$, $i = 1, \ldots, n$, such that

$$\{\mathfrak{p}\} = \bigcap_{i=1}^{n} P_{F_i|E_i}^K(\sigma_i).$$

But if $\mathfrak{p}$ is contained in $P_{K|\mathbb{Q}}^K(1)$, then this equality contradicts lemma 2.9. If $K$ is an arbitrary number field, then let $N$ be the Galois closure of $K|\mathbb{Q}$. Since

$$\varphi_{N|K} : (\mathcal{P}_N, \mathcal{T}_N) \longrightarrow (\mathcal{P}_K, \mathcal{T}_K)$$

is a continuous surjective map and $[N : K]$ is finite, $\mathcal{T}_N$ would be discrete if $\mathcal{T}_K$ is discrete.

**Proposition 3.3**
(i) Let $P_{F|E}^K(\sigma) \in \mathcal{B}_K$. Then

$$\mathcal{P}_K \setminus \left( P_{F|E}^K(\sigma) \cup R^K(F|E) \right) = \dot{\bigcup_{\langle\langle\tau\rangle\rangle \neq \langle\langle\sigma\rangle\rangle}} P_{F|E}^K(\tau),$$

and so $\varphi_{K|E}^{-1} P_{F|E}(\sigma) \cup \varphi_{K|E}^{-1} R(F|E)$ is a closed set.

(ii) Every finite subset of $\mathcal{P}_K$ is intersection of countable many elements of $\mathcal{C}_K$.

**Proof:** Assertion (i) follows from the equation

$$\mathcal{P}_K = \dot{\bigcup_{\langle\langle\tau\rangle\rangle}} \varphi_{K|E}^{-1} P_{F|E}(\tau) \,\dot{\cup}\, \varphi_{K|E}^{-1} R(F|E).$$

Let $T$ be a finite subset of $\mathcal{P}_K$. By theorem 1.3(i) we know that $T = D(K^T|K)$. Using lemma 1.2(iii), we have

$$T = D(K^T|K) = \bigcap_{j \in J} D(L_j|K),$$

where $L_j|K$ runs through the finite Galois extensions inside $K^T|K$. Thus $T$ is the intersection of countable many elements of $\mathcal{C}_K$, i.e. we proved (ii). $\square$

**Remark:** The set $P_{F|E}^K(\sigma) \cup R^K(F|E)$ is not necessarily the closure of $P_{F|E}^K(\sigma)$, since there may be isolated points in the set $R^K(F|E)$, see proposition 3.11, or there may exist subextensions of $F|E$ in which elements of $R^K(F|E)$ are unramified.

But if $K = \mathbb{Q}$, $F|\mathbb{Q}$ a Galois extension of prime degree and $\sigma \in G(F|\mathbb{Q})$, then $P_{F|\mathbb{Q}}(\sigma) \cup R(F|\mathbb{Q})$ is the closure of $P_{F|\mathbb{Q}}(\sigma)$. Suppose the contrary is true. Then there exists a prime number $p \in R(F|\mathbb{Q})$ and an open neighbourhood $U = P_{L|\mathbb{Q}}(\tau)$ of $p$, $L|\mathbb{Q}$ a finite Galois extension, such that $U$ does not meet $P_{F|\mathbb{Q}}(\sigma)$. From corollary 2.3 it follows that $F$ and $L$ are not linearly disjoint over $\mathbb{Q}$, and so $F \subseteq L$. But $p$ is unramified in $L$ and ramifies in $F$. This contradiction shows the assertion.

**Proposition 3.4**
  (i) *For every two different points $\mathfrak{p}_1$ and $\mathfrak{p}_2$ of $(\mathcal{P}_K, \mathcal{T}_K)$ there exists a both open and closed neighbourhood $W \in \mathcal{A}_K$ of $\mathfrak{p}_1$ such that $\mathfrak{p}_2 \notin W$.*
  (ii) *Let $\mathfrak{p}_1, \ldots, \mathfrak{p}_n$ be pairwise different points of $(\mathcal{P}_K, \mathcal{T}_K)$. Then there exist both open and closed neighbourhoods $U(\mathfrak{p}_i)$ of $\mathfrak{p}_i$ such that*
$$U(\mathfrak{p}_i) \cap U(\mathfrak{p}_j) = \varnothing \quad \text{for } i \neq j.$$

**Proof:** In order to prove (i) let $L|K$ be a cyclic extension of degree $m > 2$ such that $\mathfrak{p}_1$ is unramified in $L|K$ and let $\sigma \in G(L|K)$ with $\mathfrak{p}_1 \in P_{L|K}(\sigma)$. We denote the open neighbourhood $P_{L|K}(\sigma)$ of $\mathfrak{p}_1$ by $U$.

Let $N|K$ be a quadratic extension of $K$ which is unramified at all primes of $U$, completely decomposed at $R(L|K) \cup \{\mathfrak{p}_2\}$ and inert at $\mathfrak{p}_1$; if $V = P_{N|K}(\tau)$, where $\tau$ is the non-trivial element of $G(N|K)$, then $\mathfrak{p}_1 \in V$ and $\mathfrak{p}_2 \notin V$. Such an extension exists. Indeed, let

$$T = S_2 \cup S_\infty \cup R(L|K) \cup \{\mathfrak{p}_1, \mathfrak{p}_2\} \quad \text{and} \quad S = (\mathcal{P}_K \backslash U) \cup T,$$

then

$$\delta_K(S) = 1 - \frac{1}{m} > \frac{1}{2},$$

and so we can apply theorem 2.11: there exists an element $\varphi \in H^1(K_S|K, \mathbb{Z}/2\mathbb{Z})$ such that
$$\text{res}_{\mathfrak{p}}(\varphi) = 0 \in H^1(K_{\mathfrak{p}}, \mathbb{Z}/2\mathbb{Z})) \quad \text{for } \mathfrak{p} \in T \backslash \{\mathfrak{p}_1\}$$
and
$$0 \neq \text{res}_{\mathfrak{p}_1}(\varphi) \in H^1_{nr}(K_{\mathfrak{p}_1}, \mathbb{Z}/2\mathbb{Z}) \subset H^1(K_{\mathfrak{p}_1}, \mathbb{Z}/2\mathbb{Z}).$$
If $\ker \varphi = G(K_S|N)$, then $N$ is a quadratic extension of $K$ with the desired properties.

Now $W = U \cap V$ is an open neighbourhood of $\mathfrak{p}_1$ and $\mathfrak{p}_2 \notin W$. It remains to show that $W$ is closed. Let $\overline{W}$ be the closure of $W$. Using proposition 3.3(i), we get
$$\overline{U \cap V} \subseteq \overline{U} \cap \overline{V} \subseteq (U \cup R(L|K)) \cap (V \cup R(N|K)) = U \cap V,$$
and so $\overline{W} = W$. This finishes the proof of (i).

In order to prove (ii) we use induction with respect to $n$. Assume that we have found open and closed neighbourhoods $W(\mathfrak{p}_i)$ of $\mathfrak{p}_i$, $i = 1, \ldots, n-1$, which are pairwise disjoint. By (i) it follows that for every $i \in \{1, \ldots, n-1\}$ there exists an open and closed neighbourhood $W_i(\mathfrak{p}_n)$ of $\mathfrak{p}_n$ such that $\mathfrak{p}_i \notin W_i(\mathfrak{p}_n)$. Then $U(\mathfrak{p}_n) = \bigcap_{i=1}^{n-1} W_i(\mathfrak{p}_n)$ is an open and closed neighbourhood of $\mathfrak{p}_n$ such that $\mathfrak{p}_i \notin U(\mathfrak{p}_n)$ for all $i = 1, \ldots, n-1$. Now the open and closed neighbourhoods $U(\mathfrak{p}_i) = W(\mathfrak{p}_i) \setminus U(\mathfrak{p}_n)$, $i = 1, \ldots, n-1$, and $U(\mathfrak{p}_n)$ have the desired property. $\square$

Recall that a Hausdorff space $X$ is called *zero-dimensional* if every point of $X$ has a fundamental system of neighbourhoods which are both open and closed, and $X$ is called *strongly zero-dimensional* if for every closed subset $A$ of $X$ and each neighbourhood $U$ of $A$ there is an open and closed neighbourhood of $A$ contained in $U$.

**Proposition 3.5** *The space $(\mathcal{P}_K, \mathcal{T}_K)$ has the following properties: it is*
 (i) *a Hausdorff space,*
 (ii) *strongly zero-dimensional (and so totally disconnected),*
 (iii) *metrizable (and so normal and completely regular),*
 (iv) *and every point of $(\mathcal{P}_K, \mathcal{T}_K)$ has a base of neighbourhoods consisting of open and closed sets.*

**Proof:** By proposition 3.4(i) there exists for every two different points $x$ and $y$ of $\mathcal{P}_K$ an open and closed neighbourhood $W$ of $x$ such that $y \notin W$. It follows that $\mathcal{P}_K \setminus W$ is an open neighbourhood of $y$ being disjoint to $W$. Therefore $\mathcal{P}_K$ is a Hausdorff space.

Now we proved (iv). Let $\mathfrak{p} \in (\mathcal{P}_K, \mathcal{T}_K)$ and let
$$U = \bigcap_{i=1}^{n} P^K_{F_i|E_i}(\sigma_i) \in \mathcal{A}_K$$
be an open neighbourhood of $\mathfrak{p}$. We have to find and open and closed neighbourhood of $\mathfrak{p}$ being contained in $U$. Obviously we may assume that $U = P^K_{F|E}(\sigma)$. By proposition 3.4(ii) there exist open und closed, pairwise disjoint neighbourhoods $U(\mathfrak{p}_i)$ of $\mathfrak{p}_i$ $i = 0, \ldots, n$, where $\{\mathfrak{p}_1 \ldots, \mathfrak{p}_n\} = R^K(F|E)$ and $\mathfrak{p}_0 = \mathfrak{p}$. Then

$$U_R = \dot{\bigcup_{i=1}^{n}} U(\mathfrak{p}_i)$$

is an open and closed neighbourhood of $R^K(F|E)$ not containing $\mathfrak{p}$. Let

$$V = P^K_{F|E}(\sigma) \backslash U_R,$$

then $V$ is open and contains $\mathfrak{p}$. But $V$ is also closed, since we get for the closure $\overline{V}$ of $V$, using proposition 3.3(i),

$$\begin{aligned}
\overline{P^K_{F|E}(\sigma) \backslash U_R} &= \overline{P^K_{F|E}(\sigma) \cap (\mathcal{P}_K \backslash U_R)} \\
&\subseteq \overline{P^K_{F|E}(\sigma)} \cap (\mathcal{P}_K \backslash U_R) \\
&\subseteq (P^K_{F|E}(\sigma) \cup R^K(F|E)) \cap (\mathcal{P}_K \backslash U_R) \\
&= P^K_{F|E}(\sigma) \backslash U_R.
\end{aligned}$$

This finishes the proof of (iv). The other assertions follow from [2] IX.6 exercise 2(b) since the considered space has a countable base. □

**Remarks:**
1. The space $(\mathcal{P}_K, \mathcal{T}_K)$ is not compact. Otherwise $(\mathcal{P}_\mathbb{Q}, \mathcal{T}_\mathbb{Q})$ would be compact. But the following is true:

*Every compact subset $A$ of $\mathcal{P}_\mathbb{Q}$ has no interior point.*

Suppose the contrary, i.e. there is an open subset $U$ of $\mathcal{P}_\mathbb{Q}$ contained in $A$. Then, by the remark following definition 3.1, we may assume that $U$ is a Čebotarev set: $U = P_{F|\mathbb{Q}}(\sigma)$. Let $L_{p_1}, \ldots, L_{p_n}$ be the finitely many quadratic fields with prime number discriminant contained in $F$; here we use the notation of proposition 2.6. Using this proposition, we get

$$\bigcup_{p \neq p_1, \ldots, p_n} (P_{L_p|\mathbb{Q}}(1) \cap A) = A,$$

and so there finitely many prime numbers $q_1, \ldots, q_m \notin \{p_1, \ldots, p_n\}$ such that

$$\bigcup_{i=1}^{m} (P_{L_{q_i}|\mathbb{Q}}(1) \cap A) = A, \quad \text{i.e.} \quad A \subseteq \bigcup_{i=1}^{m} P_{L_{q_i}|\mathbb{Q}}(1).$$

Since

$$\left( P_{L_{q_1}|\mathbb{Q}}(1) \cup \cdots \cup P_{L_{q_m}|\mathbb{Q}}(1) \right) \cap \left( P_{L_{q_1}|\mathbb{Q}}(-1) \cap \cdots \cap P_{L_{q_m}|\mathbb{Q}}(-1) \right) = \emptyset,$$

and, by the propositions 2.2 and 2.3,

$$A \cap \bigcap_{i=1}^{m} P_{L_{q_i}|\mathbb{Q}}(-1) \supseteq P_{F|\mathbb{Q}}(\sigma) \cap \bigcap_{i=1}^{m} P_{L_{q_i}|\mathbb{Q}}(-1) \neq \emptyset,$$

we get a contradiction.

2. If the set $\mathcal{P}'_K = \mathcal{P}_K \cup \{(0)\}$ of all prime ideals of the ring of integers $\mathcal{O}_K$ of $K$ is equipped with the sum of the topologies $\mathcal{T}_K$ of $\mathcal{P}_K$ and the discrete topology on $\{(0)\}$, then the identity map

$$(\mathcal{P}'_K, \mathcal{T}'_K) \xrightarrow{id} \operatorname{Spec}_{Zar} \mathcal{O}_K$$

is continuous, since the non-trivial closed subsets of $\operatorname{Spec} \mathcal{O}_K$ with respect to the Zariski topology are the finite sets not containing $(0)$ and $(\mathcal{P}_K, \mathcal{T}_K)$ is a Hausdorff space.

**Definition 3.6** *An open set $U$ of $(\mathcal{P}_K, \mathcal{T}_K)$ is called* **finitely presented**, *if it is a finite union of elements of $\mathcal{A}_K$, i.e.*

$$U = \bigcup_{i=1}^{n} \bigcap_{j=1}^{n_i} P^K_{F_{ij}|E_{ij}}(\sigma_{ij}).$$

Examinating the proofs of the propositions 3.5(iv) and 3.4 we get

**Lemma 3.7** *Every point of $(\mathcal{P}_K, \mathcal{T}_K)$ has a fundamental system of neighbourhoods which are both open and closed and finitely presented, and we can separate finitely many points by disjoint open and closed and finitely presented sets.*

**Lemma 3.8** *If $(X, \mathcal{T})$ is a strongly zero-dimensional topological space, then for every finite open covering $\bigcup_{i=1}^{n} U_i = X$ of $X$ there exist pairwise disjoint, open and closed sets $V_i \subseteq U_i$, $i = 1, \ldots, n$, such that*

$$\dot{\bigcup}_{i=1}^{n} V_i = X.$$

**Proof:** Assume that we have found open and closed set $W_i \subseteq U_i$, $i = 1, \ldots, m$, with $m \leq n$ such that

$$\bigcup_{i=1}^{m} W_i \cup \bigcup_{i=m+1}^{n} U_i = X.$$

If $A = X \setminus (\bigcup_{i=1}^{m} W_i \cup \bigcup_{i=m+2}^{n} U_i)$, then $A$ is closed and $A \subseteq U_{m+1}$. Since $X$ is strongly zero-dimensional, there exists an open and closed set $W_{m+1}$ with $A \subseteq W_{m+1} \subseteq U_{m+1}$. It follows that

$$\bigcup_{i=1}^{m+1} W_i \cup \bigcup_{i=m+2}^{n} U_i = X,$$

and finally $\bigcup_{i=1}^{n} W_i = X$. Now the sets

$$V_i = W_i \setminus \bigcup_{j<i} W_j \subseteq U_i$$

are open and closed, pairwise disjoint and they form a covering of $X$. □

From the lemma above and proposition 3.5 we obtain the

**Corollary 3.9** *If $\bigcup_{i=1}^{n} U_i = \mathcal{P}_K$ is a finite open covering of the space $(\mathcal{P}_K, \mathcal{T}_K)$, then there exist pairwise disjoint, open and closed sets $V_i \subseteq U_i$, $i = 1, \ldots, n$, such that*

$$\dot{\bigcup}_{i=1}^{n} V_i = \mathcal{P}_K.$$

**Proposition 3.10**
  (i) *Let $\mathfrak{p} \in (\mathcal{P}_K, \mathcal{T}_K)$ be a prime ideal of $K$ such that $p = \mathfrak{p} \cap \mathbb{Q}$ is completely decomposed in $K$. Then every open neighbourhood of $\mathfrak{p}$ has positive density.*
  (ii) *Let $\mathfrak{p} \in (\mathcal{P}_K, \mathcal{T}_K)$ be a prime ideal of $K$ such that $p = \mathfrak{p} \cap \mathbb{Q}$ is unramified in $K$. Then every open neighbourhood of $\mathfrak{p}$ has infinitely many points.*

**Proof:** Let $\mathfrak{p} \in (\mathcal{P}_K, \mathcal{T}_K)$ such that $p = \mathfrak{p} \cap \mathbb{Q}$ is completely decomposed in $K$ and let $U$ be an open neighbourhood of $\mathfrak{p}$. The prime number $p$ is also completely decomposed in the normal closure $N$ of $K|\mathbb{Q}$. If $\mathfrak{P}$ is an extension of $\mathfrak{p}$ to $N$, then $V = \varphi_{N|K}^{-1}(U)$ is an open neighbourhood of $\mathfrak{P}$. Since every open neighbourhood of a point of $(\mathcal{P}_N, \mathcal{T}_N)$ contains a set which is a finite intersection of sets of $\mathcal{B}_N$, it follows from lemma 2.9 that $V$ has positive density, and so $U$ has. This proves assertion (i) and (ii) follows from lemma 2.10. □

Recall that a point $x$ of a topological space $X$ is called *isolated* if $\{x\}$ is an open set in $X$.

If $G(F|E)$ is the Galois group of a finite Galois extension $F|E$ and $\mathfrak{P}$ a prime of $F$, then we denote the decomposition group and the inertia subgroup of $G(F|E)$ with respect to $\mathfrak{P}$ by $G_\mathfrak{P} = G_\mathfrak{P}(F|E)$ and $T_\mathfrak{P} = T_\mathfrak{P}(F|E)$, respectively. If $\ell$ is a prime number, then $G(\ell)$ is a $\ell$-Sylow group of a group $G$.

**Proposition 3.11** *Let $K|\mathbb{Q}$ be a finite extension and let $\mathfrak{p} \in \varphi_{K|\mathbb{Q}}^{-1}(R(K|\mathbb{Q}))$.*
(i) *Assume that $K|\mathbb{Q}$ is normal and that $G_\mathfrak{p}(K|\mathbb{Q})$ has the following property:*

  *there exists a prime number $\ell$ such that $G_\mathfrak{p}(\ell)$ is not cyclic and the quotient $G_\mathfrak{p}(\ell)/T_\mathfrak{p}(\ell)$ is non-trivial.*

  *Then $\mathfrak{p}$ is an isolated point of $(\mathcal{P}_K, \mathcal{T}_K)$.*
(ii) *For every prime ideal $\mathfrak{P}|\mathfrak{p}$ of the normal closure $N$ of $K|\mathbb{Q}$ there exists a finite Galois extension $L|N$ such that $\mathfrak{P}$ and all $G(N|\mathbb{Q})$-conjugates of $\mathfrak{P}$ are inert in $L|N$ and their unique extensions to $L$ are isolated in $(\mathcal{P}_L, \mathcal{T}_L)$.*

**Proof:** Let $K_0 \subseteq K$ be the fixed field of $[G_\mathfrak{p}(\ell), G_\mathfrak{p}(\ell)]$. From our assumptions it follows that $K_0$ has subfields $E_i$, $i = 0, 1, 2$, such that $K_0 = E_1 E_2$, $E_0 = E_1 \cap E_2$ and $G(K_0|E_0) \cong \mathbb{Z}/\ell\mathbb{Z} \times \mathbb{Z}/\ell\mathbb{Z}$, and $\mathfrak{p} \cap E_1$ is inert and $\mathfrak{p} \cap E_2$ is ramified in $K_0$. Let $E_3$ be any extension of $E_0$ in $K_0$ of degree $\ell$ different to $E_1$ and $E_2$:

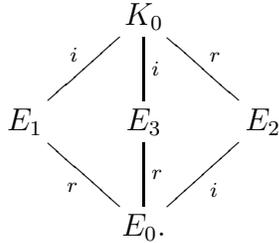

The letters $i$ and $r$ indicate whether $\mathfrak{p} \cap E_0$ resp. its unique extensions to the fields $E_i$, $i = 1, 2, 3$, are inert or ramify in the considered extensions. Now we consider the open set

$$U = P_{K_0|E_1}^{K_0}(\sigma) \cap P_{K_0|E_3}^{K_0}(\tau) = \varphi_{K_0|E_1}^{-1} P_{K_0|E_1}(\sigma) \cap \varphi_{K_0|E_3}^{-1} P_{K_0|E_3}(\tau)$$

of $(\mathcal{P}_{K_0}, \mathcal{T}_{K_0})$, where

$$\sigma = \Big(\frac{K_0|E_1}{\mathfrak{p} \cap K_0}\Big), \qquad \tau = \Big(\frac{K_0|E_3}{\mathfrak{p} \cap K_0}\Big).$$

Observe that $\sigma \neq 1 \neq \tau$ and $\mathfrak{p}_0 = \mathfrak{p} \cap K_0 \in U$.

Let $\mathfrak{p}'$ be a prime ideal contained in $U$. Since $K_0|E_0$ is not cyclic and $\mathfrak{p}' \cap E_1$ is inert in $K_0|E_1$, $\mathfrak{p}' \cap E_0$ is completely decomposed or ramifies in $E_1|E_0$. In the first case its extensions to $E_3$ would also be completely decomposed in $K_0|E_3$, and so $\mathfrak{p}'$ can not be contained in $\varphi_{K_0|E_3}^{-1} P_{K_0|E_3}(\tau)$. It follows that $U \subseteq \varphi_{K_0|E_0}^{-1} R(K_0|E_0)$, and so $U$ is finite. Therefore $\{\mathfrak{p}_0\} \subseteq U$ is also open (the finite set $U \setminus \{\mathfrak{p}_0\}$ is

closed as $(\mathcal{P}_{K_0}, \mathcal{T}_{K_0})$ is a Hausdorff space). Therefore $\mathfrak{p}_0$ is an isolated point of $(\mathcal{P}_{K_0}, \mathcal{T}_{K_0})$, and so $\mathfrak{p} = \varphi_{K|K_0}^{-1}(\mathfrak{p}_0)$ is an isolated point of $(\mathcal{P}_K, \mathcal{T}_K)$. This proves assertion (i).

In order to prove (ii) let $\mathfrak{P}$ be a prime ideal contained in $\varphi_{N|\mathbb{Q}}^{-1}(R(N|\mathbb{Q}))$ and let $\ell$ be any prime number dividing the order of the inertia subgroup $T_\mathfrak{P}$ of $G_\mathfrak{P} = G_\mathfrak{P}(N|\mathbb{Q})$. Let $L_0|\mathbb{Q}$ be a cyclic extension of $\ell$-power degree such that $\mathfrak{P} \cap \mathbb{Q}$ is inert in $L_0|\mathbb{Q}$ and $L_0 \subsetneq N$. Let $L = NL_0$. Then all $G(N|\mathbb{Q})$-conjugates of $\mathfrak{P}$ are inert in $L|N$ and $G_{\mathfrak{P}_L}(L|\mathbb{Q})$ fulfills the condition of (i), where $\mathfrak{P}_L$ denotes the unique extension of $\mathfrak{P}$ to $L$. It follows that $\mathfrak{P}_L$ is isolated in $(\mathcal{P}_L, \mathcal{T}_L)$. $\square$

**Definition 3.12** *Let $K$ be a number field and $N$ the normal closure of $K|\mathbb{Q}$. A point $\mathfrak{p} \in (\mathcal{P}_K, \mathcal{T}_K)$ is called **potentiell isolated** if for every $\mathfrak{P}|\mathfrak{p}$ of $N$ there exists a finite Galois extension $L|N$ such that*

(i) *all $G(N|\mathbb{Q})$-conjugates of $\mathfrak{P}$ are unramified in $L|N$,*

(ii) *all points of $\varphi_{L|N}^{-1}(\mathfrak{P})$ are isolated in $(\mathcal{P}_L, \mathcal{T}_L)$.*

*We denote the set of all isolated points and the set of all potentiell isolated points of $(\mathcal{P}_K, \mathcal{T}_K)$ by $(\mathcal{P}_K)_{iso}$ and $(\mathcal{P}_K)_{p.iso}$, respectively.*

Without condition (i) in the definition above, i.e. $\varphi_{N|\mathbb{Q}}^{-1}(\mathfrak{P} \cap \mathbb{Q}) \subseteq U(L|N)$, all points of $\mathcal{P}_K$ would be potentiell isolated, since for every $\mathfrak{p} \in \mathcal{P}_K$ there exists a finite Galois extension $K'|K$ in which $\mathfrak{p}$ ramifies, and we can apply proposition 3.11(ii) to the field $K'$. Furthermore we would like to mention (although it is completely trivial) that

$$\mathcal{P}_\mathbb{Q} \text{ has no isolated points,}$$

since every open set of $\mathcal{P}_\mathbb{Q}$ has positive density. The following proposition considers the general case.

**Theorem 3.13** *Let $K$ be a number field. Then the following is true:*

(i) $$(\mathcal{P}_K)_{iso} \subseteq \varphi_{K|\mathbb{Q}}^{-1}(R(K|\mathbb{Q})) = (\mathcal{P}_K)_{p.iso},$$

(ii) $\varphi_{K|\mathbb{Q}}^{-1}(U(K|\mathbb{Q})) \subseteq \{\mathfrak{p} \in \mathcal{P}_K \mid$ *every open neighbourhood of $\mathfrak{p}$ has infinitely many points* $\}$,

(iii) $\varphi_{K|\mathbb{Q}}^{-1}(D(K|\mathbb{Q})) \subseteq \{\mathfrak{p} \in \mathcal{P}_K \mid$ *every open neighbourhood of $\mathfrak{p}$ has positive density* $\}$.

*If $K|\mathbb{Q}$ is a Galois extension, then we have equality in (iii).*

**Proof:** Let $N$ be the normal closure over $K$ over $\mathbb{Q}$. The inclusion

$$\varphi_{K|\mathbb{Q}}^{-1}(R(K|\mathbb{Q})) \subseteq (\mathcal{P}_K)_{p.iso}$$

is just proposition 3.11(ii). Suppose that $\mathfrak{p} \in (\mathcal{P}_K)_{p.iso}$ is not contained in $\varphi_{K|\mathbb{Q}}^{-1}(R(K|\mathbb{Q}))$. Then the extensions $\mathfrak{P}$ of $\mathfrak{p}$ to $N$ are contained in $\varphi_{N|\mathbb{Q}}^{-1}(U(N|\mathbb{Q}))$. Let $\mathfrak{P}_0$ be one of these extensions and let $L|N$ be a finite Galois extension such that all $G(N|\mathbb{Q})$-conjugates of $\mathfrak{P}_0$ are unramified in $L|N$ and all points $\mathfrak{P}_{0L} \in \varphi_{L|N}^{-1}(\mathfrak{P}_0)$ are isolated in $(\mathcal{P}_L, \mathcal{T}_L)$. Then $\mathfrak{P}_{0L} \in \varphi_{L|\mathbb{Q}}^{-1}(U(L|\mathbb{Q}))$. This contradicts proposition 3.10(ii) and therefore we proved (i).

Assertions (ii) (and so the inclusion in (i)) and the inclusion (iii) follow from proposition 3.10(ii) and (i), respectively.

Now we show that for every point $\mathfrak{P} \in \varphi_{N|\mathbb{Q}}^{-1}(U(N|\mathbb{Q})) \setminus \varphi_{N|\mathbb{Q}}^{-1}(D(N|\mathbb{Q}))$ there exists an open neighbourhood of density equal to 0. Indeed, let $N_0 \subset N$ be its decomposition field and observe that by assumption $N \neq N_0$. Therefore $\tau = \left(\frac{N|N_0}{\mathfrak{P}}\right) \in G(N|N_0)$ is not equal to 1. Obviously, $\mathfrak{P} \in \varphi_{N|N_0}^{-1} P_{N|N_0}(\tau)$ and this open set has density equal to 0 since every prime ideal of $P_{N|N_0}(\tau)$ is inert in the extension $N|N_0$. So we get

$$\varphi_{N|\mathbb{Q}}^{-1}(D(N|\mathbb{Q})) = \left\{\mathfrak{P} \in \mathcal{P}_N \mid \begin{array}{l}\text{every open neighbourhood of }\mathfrak{P}\\ \text{has positive density}\end{array}\right\}$$

showing also the equality stated in (iii). $\square$

**Remark 1:** The inclusion in (ii) may be strict (even if $K|\mathbb{Q}$ is a Galois extension), i.e. there may exist ramified primes having only infinite open neighbourhoods, or with other words, it is possible that there are ramified points which are not isolated. An example of this situation is the following (see also proposition 3.11(i) where we consider the opposite situation):

Let $K|\mathbb{Q}$ be a cyclic extension of prime degree ramified at $p$. Suppose that the unique extension $\mathfrak{p}$ of $p$ to $K$ is an isolated point. Then there exist finitely many Galois extensions $L|K$ and $E|\mathbb{Q}$ and elements $\sigma \in G(L|K)$, $\tau \in G(E|\mathbb{Q})$, such that

$$\{\mathfrak{p}\} = \bigcap_{L|K} P_{L|K}(\sigma) \cap \bigcap_{E|\mathbb{Q}} \varphi_{K|\mathbb{Q}}^{-1} P_{E|\mathbb{Q}}(\tau).$$

Using corollary 2.2, it follows that

$$\{\mathfrak{p}\} = P_{L|K}(\sigma) \cap \varphi_{K|\mathbb{Q}}^{-1} P_{E|\mathbb{Q}}(\tau)$$

for one extension $L|K$ and one extension $E|\mathbb{Q}$. Furthermore, we may assume that $E \subseteq L$. This follows from proposition 2.1 with $N = LE$, since the extensions of $\mathfrak{p}$ to $L$ are unramified in $LE$ (as $\mathfrak{p} \in \varphi_{K|\mathbb{Q}}^{-1} P_{E|\mathbb{Q}}(\tau) \subseteq \varphi_{K|\mathbb{Q}}^{-1} U(E|\mathbb{Q})$). So we have the following diagram of fields:

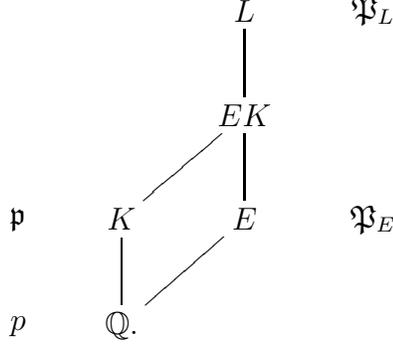

We claim that without lost of generality $\sigma_{|E} = \tau$.

Since $\mathfrak{p} \in P_{L|K}(\sigma) \cap \varphi_{K|\mathbb{Q}}^{-1} P_{E|\mathbb{Q}}(\tau)$, there exist primes $\mathfrak{P}_L | \mathfrak{p}$ of $L$ and $\mathfrak{P}_E | p$ of $E$ such that
$$\sigma = \left(\frac{L|K}{\mathfrak{P}_L}\right) \quad \text{and} \quad \tau = \left(\frac{E|\mathbb{Q}}{\mathfrak{P}_E}\right).$$

Since the residue degree $f(\mathfrak{p}|p) = 1$, we have
$$\left(\frac{L|\mathbb{Q}}{\mathfrak{P}_L}\right) = \left(\frac{L|\mathbb{Q}}{\mathfrak{P}_L}\right)^{f(\mathfrak{p}|p)} = \left(\frac{L|K}{\mathfrak{P}_L}\right) = \sigma$$

and so
$$\sigma_{|E} = \left(\frac{L|\mathbb{Q}}{\mathfrak{P}_L}\right)_{|E} = \left(\frac{E|\mathbb{Q}}{\mathfrak{P}_L \cap E}\right) = \left(\frac{E|\mathbb{Q}}{\rho \mathfrak{P}_E}\right) = \tau^\rho$$

for some $\rho \in G(E|\mathbb{Q})$. This shows the claim.

Since $p$ ramifies in $K|\mathbb{Q}$, the subset
$$V = P_{L|K}(\sigma) \cap \varphi_{K|\mathbb{Q}}^{-1}\left(P_{E|\mathbb{Q}}(\tau) \cap P_{K|\mathbb{Q}}(1)\right)$$

of $P_{L|K}(\sigma) \cap \varphi_{K|\mathbb{Q}}^{-1} P_{E|\mathbb{Q}}(\tau) = \{\mathfrak{p}\}$ does not contain $\mathfrak{p}$, and so $V = \varnothing$. But this is a contradiction, since for a prime $\mathfrak{p}'$ in $P_{L|K}(\sigma) \cap \varphi_{K|\mathbb{Q}}^{-1} P_{K|\mathbb{Q}}(1)$ (this set has positive density) it follows that there is a prime $\mathfrak{P}'_L | \mathfrak{p}'$ of $L$ such that $\left(\frac{L|K}{\mathfrak{P}'_L}\right) = \sigma$, and so
$$\left(\frac{EK|K}{\mathfrak{P}'_L \cap EK}\right) = \left(\frac{L|K}{\mathfrak{P}'_L}\right)_{|EK} = \sigma_{|EK}$$

from which follows that
$$\left(\frac{E|\mathbb{Q}}{\mathfrak{P}'_L \cap E}\right) = \left(\frac{E|\mathbb{Q}}{\mathfrak{P}'_L \cap E}\right)^{f(\mathfrak{p}'|p)} = \left(\frac{EK|K}{\mathfrak{P}'_L \cap EK}\right)_{|E} = \sigma_{|E} = \tau.$$

Thus $\mathfrak{p}' \in \varphi_{K|\mathbb{Q}}^{-1}(P_{E|\mathbb{Q}}(\tau))$, and so $\mathfrak{p}' \in V$. It follows that $\mathfrak{p}$ is not an isolated point.

**Remark 2:** We will show that

$$\mathcal{P}_K \text{ has infinite compact subsets.}$$

Indeed, let $x \in \mathcal{P}_K$ be a point not contained in $(\mathcal{P}_K)_{p.iso}$ and let

$$\mathcal{V} = \{V_i \mid i \in \mathbb{N}, V_i \text{ both open and closed}, x \in V_i\}$$

be a base of open and closed neighbourhoods of $x$ such that

$$V_{i+1} \subseteq V_i, \ V_{i+1} \neq V_i, \quad i \in \mathbb{N}.$$

Such a base exists, since the topology $\mathcal{T}_K$ is countable and every open neighbourhood of $x$ contains infinitely many points by theorem 3.13 and we can separate two different points by both open and closed sets. Let $x_i \in V_i \backslash V_{i+1}$ and

$$A = \{x\} \cup \{x_1, x_2, \ldots\}.$$

Then $A$ is infinite and compact.

**Proposition 3.14** *Let $K|\mathbb{Q}$ be a number field. Then there exists a finite Galois extension $L|K$ such that*

$$\varphi_{L|\mathbb{Q}}^{-1}(R(L|\mathbb{Q})) = (\mathcal{P}_L)_{iso}.$$

**Proof:** Without lost of generality we may assume that $K|\mathbb{Q}$ is normal. Let $\mathfrak{p} \in \varphi_{K|\mathbb{Q}}^{-1}(R(K|\mathbb{Q}))$, $\ell$ a prime number dividing the order of the inertia group $T_{\mathfrak{p}}(K|\mathbb{Q})$ and $K' = K(\mu_{2\ell})$. Furthermore let $q_1, \ldots, q_t$, $t > [K' : \mathbb{Q}]/2 + 1$, be auxiliary pairwise different prime numbers not above $\ell$,

$$S = \{\mathfrak{q} \mid \mathfrak{q} = \mathfrak{p} \text{ or } \mathfrak{q} | \ell \cdot q_1 \cdots q_t\} \cup S_\infty$$

and $K'_S(\ell)$ the maximal $\ell$-extension of $K'$ unramified outside $S$. From [4] (10.7.8), (10.7.9) and (10.6.4) it easily follows that the local groups $G(K'_{\mathfrak{q}}(\ell)|K'_{\mathfrak{q}})$ of the maximal $\ell$-extension $K'_{\mathfrak{q}}(\ell)$ of $K'_{\mathfrak{q}}$, $\mathfrak{q} \in S \backslash S_\infty$, are contained in $G(K'_S(\ell)|K')$ (if the group $G(K'_S(\ell)|K')$ would be degenerated, then its $\mathbb{Z}_\ell$-rank has to be equal to the number of complex places of $K'$ plus 1). Since these groups are not (pro-)cyclic, there exists a finite Galois extension $L^{\mathfrak{p}}$ of $K'$ such that the decomposition groups of all primes of $S(L^{\mathfrak{p}})$ fulfill the assumptions of proposition 3.11(i). Now the finite Galois extension

$$L = \prod_{\mathfrak{p} \in R^K(K|\mathbb{Q})} L^{\mathfrak{p}}$$

of $K$ has the property that all points of $R^L(L|\mathbb{Q})$ are isolated in $(\mathcal{P}_L, \mathcal{T}_L)$. □

Now we consider the inverse limit topology $\mathcal{T}_{\bar{K}}$ on the space

$$\mathcal{P}_{\bar{K}} = \varprojlim_{K'|K} \mathcal{P}_{K'}$$

where the inverse limit is taken over all spaces $(\mathcal{P}_{K'}, \mathcal{T}_{K'})$, where $K'|K$ runs through the finite subextensions contained in a fixed algebraic closure $\bar{K}$ of $K$ and the transition maps are the continuous maps $\varphi_{K'|K}$.

**Proposition 3.15** *The space $(\mathcal{P}_{\bar{K}}, \mathcal{T}_{\bar{K}})$ is a locally compact, totally disconnected Hausdorff space which is countable at infinity. The Alexandroff-compactification $\mathcal{P}_{\bar{K}} \cup \{\omega\}$ of $(\mathcal{P}_{\bar{K}}, \mathcal{T}_{\bar{K}})$ is a profinite space, i.e. it is compact and totally disconnected.*

**Proof:** Since the spaces $(\mathcal{P}_{K'}, \mathcal{T}_{K'})$, $K'|K$ finite, are totally disconnected Hausdorff spaces, the same is true for their inverse limit.

Let $x_{K'}$ be a point of $(\mathcal{P}_{K'}, \mathcal{T}_{K'})$, $K'|K$ finite. Then there exists a finite extension of $K'$ such that $x_{K'}$ ramifies in this extension, and so, by proposition 3.14, there exists a finite extension $K''|K'$ such that the set $\varphi_{K''|K'}^{-1}(x_{K'})$ consists of isolated points of $(\mathcal{P}_{K''}, \mathcal{T}_{K''})$. If we denote the canonical (continuous) map from the inverse limit to $\mathcal{P}_{K'}$ by

$$\varphi_{K'} : (\mathcal{P}_{\bar{K}}, \mathcal{T}_{\bar{K}}) \longrightarrow (\mathcal{P}_{K'}, \mathcal{T}_{K'}),$$

then it follows that the set

$$\varphi_{K'}^{-1}(x_{K'}) = \varphi_{K''}^{-1} \varphi_{K''|K'}^{-1}(x_{K'})$$

is open and closed in $(\mathcal{P}_{\bar{K}}, \mathcal{T}_{\bar{K}})$. Since $\varphi_{K'''|K'}^{-1}(x_{K'})$ is finite for every finite extension $K'''|K''$, it follows that $\varphi_{K'}^{-1}(x_{K'})$ is compact and $\varphi_{K'}$ is surjective, see [2] I.9.6, cor.1. It follows that $\mathcal{P}_{\bar{K}}$ is the union of the open compact sets $\varphi_K^{-1}(x_K)$, $x_K \in \mathcal{P}_K$. Therefore $(\mathcal{P}_{\bar{K}}, \mathcal{T}_{\bar{K}})$ is locally compact and countable at infinity because $\mathcal{P}_K$ is countable.

Finally we show that the Alexandroff-compactification $\mathcal{P}_{\bar{K}} \cup \{\omega\}$ is totally disconnected. Indeed, the sets $(\mathcal{P}_{\bar{K}} \setminus \varphi_K^{-1}(x)) \cup \{\omega\}$, $x \in \mathcal{P}_K$, are both open and closed neighbourhoods of $\omega$, and

$$\bigcap_{x \in \mathcal{P}_K} \left( (\mathcal{P}_{\bar{K}} \setminus \varphi_K^{-1}(x)) \cup \{\omega\} \right) = \{\omega\}.$$

Therefore the connected component of $\omega$ is equal to $\{\omega\}$. $\square$

# 4 Uniformity

In this section we consider uniformities on $\mathcal{P}_K$ which induce the topology $\mathcal{T}_K$. First we recall some facts concerning uniform structures on a normal topological space $(X, \mathcal{T})$.

*1. The uniformity of finite partitions by open and closed subsets of $X$.*

A partition of $X$ is a family $(U_i)_{i \in I}$ of pairwise disjoint subsets such that their union is the whole space. Let

$$Part^{oc} = \{(V_1, \ldots, V_n) \mid n \in \mathbb{N},\ V_i \subseteq (X, \mathcal{T}) \text{ open and closed},\ \dot{\bigcup_{i=1}^n} V_i = X\}$$

be the set of finite partitions of $X$ by open and closed subsets of $(X, \mathcal{T})$ and let

$$\mathfrak{V}^{oc} = \{\, V_Q = \dot{\bigcup_{i=1}^n}(V_i \times V_i) \subseteq X \times X \mid n \in \mathbb{N},\ Q = (V_1, \ldots, V_n) \in Part^{oc}\}.$$

Obviously, $\mathfrak{V}^{oc}$ is a base for the uniform structure

$$\mathfrak{U}^{oc} = \{U \subseteq X \times X \mid U \text{ contains a subset } V \in \mathfrak{V}^{oc}\}$$

on $X$. We denote the completion of $(X, \mathfrak{U}^{oc})$ by $(\hat{X}, \hat{\mathfrak{U}}^{oc})$.

*2. The uniformity of finite open coverings on $X$.*

Let

$$Cov^o = \{(U_1, \ldots, U_n) \mid n \in \mathbb{N},\ U_i \in \mathcal{T},\ \bigcup_{i=1}^n U_i = X\}$$

be the set of finite coverings of $X$ by open subsets of $(X, \mathcal{T})$ and let

$$\mathfrak{V}^o = \{\, V_Q = \bigcup_{i=1}^n (U_i \times U_i) \subseteq X \times X \mid n \in \mathbb{N},\ Q = (U_1, \ldots, U_n) \in Cov^o\}.$$

This is a base for the uniform structure

$$\mathfrak{U}^o = \{U \subseteq X \times X \mid U \text{ contains a subset } V \in \mathfrak{V}^o\}$$

on $X$, the so-called uniformity of finite open coverings.

3. The Stone-Čech compactification of $X$.

Let $\beta X$ be the Stone-Čech compactification of $(X, \mathcal{T})$. Then $\beta X$ is the completion of $X$ with respect to the coarsest uniformity $\mathfrak{U}^{S\check{C}}$ on $X$ for which all continuous mappings of $X$ into $[0,1]$ are uniformly continuous (see [2] IX.1 exercise 7).

Concerning these three uniformities on $X$ we have the following

**Proposition 4.1** *Let $(X, \mathcal{T})$ be a strongly zero-dimensional Hausdorff space. Then the uniform structures $\mathfrak{U}^{oc}$, $\mathfrak{U}^o$ and $\mathfrak{U}^{S\check{C}}$ on $X$ are equal, now denoted by $\mathfrak{U}$. The topology induced by $\mathfrak{U}$ on $X$ is equal to $\mathcal{T}$.*

**Proof:** Since the space $(X, \mathcal{T})$ is normal, the uniformity $\mathfrak{U}^o$ is equal to the uniformity $\mathfrak{U}_K^{S\check{C}}$ and $\mathfrak{U}^o = \mathfrak{U}^{S\check{C}}$ induces the topology $\mathcal{T}$ on $X$, see [2] IX.4 exercise 17(b).

By definition $\mathfrak{U}^{oc}$ is coarser than $\mathfrak{U}^o$ and, using lemma 3.8, there exists for every $V_{(U_1,\ldots,U_n)} \in \mathfrak{V}^o$ an element $V_{(V_1,\ldots,V_n)} \in \mathfrak{V}^{oc}$ with $V_i \subseteq U_i$, $i = 1, \ldots, n$, i.e. $V_{(V_1,\ldots,V_n)} \subseteq V_{(U_1,\ldots,U_n)}$. Thus $\mathfrak{U}^{oc}$ is finer than $\mathfrak{U}^o$, and so they are equal. $\square$

**Proposition 4.2** *Let $X$ be a strongly zero-dimensional Hausdorff space. Then the completion $(\hat{X}, \hat{\mathfrak{U}})$ of $X$ equipped with the uniformity $\mathfrak{U} = \mathfrak{U}^{oc}$ is a profinite space, i.e. it is compact and totally disconnected.*

**Proof:** From proposition 4.1 it follows that $(\hat{X}, \hat{\mathfrak{U}}) = \beta(X, \mathcal{T})$ and the compact space $\beta(X, \mathcal{T})$ is totally disconnected, see [2] IX.6 exercise 1(b). $\square$

**Proposition 4.3** *Let $X$ be a strongly zero-dimensional Hausdorff space and let*

$$i : (X, \mathfrak{U}) \longrightarrow (\hat{X}, \hat{\mathfrak{U}})$$

*be the canonical mapping ($\mathfrak{U} = \mathfrak{U}^{oc}$) and we identify $X$ with $i(X)$. Let $\mathcal{OC}_X$ and $\mathcal{OC}_{\hat{X}}$ be the set of both open and closed subsets of $X$ and $\hat{X}$, respectively.*

(i) *The maps*

$$\mathcal{OC}_X \longrightarrow \mathcal{OC}_{\hat{X}}, \quad S \mapsto \overline{S}, \quad \text{and} \quad \mathcal{OC}_{\hat{X}} \longrightarrow \mathcal{OC}_X, \quad \mathcal{S} \mapsto \mathcal{S} \cap X$$

*are bijections, where $\overline{S}$ is the closure of $S$ in $\hat{X}$.*

(ii) *For the set of isolated points of $X$ and $\hat{X}$ we have $i(X_{iso}) = \hat{X}_{iso}$.*

**Proof:** It is clear that the second map is well-defined. Let $S \in \mathcal{OC}_X$. Since $S$ and $X \backslash S$ are closed sets of $X$, we get from $S \cup (X \backslash S) = X$ the partition

$$\overline{S} \cup \overline{X \backslash S} = \hat{X},$$

see [2] IX.4 exercise 17(c), and so $\hat{X} \backslash \overline{S} = \overline{X \backslash S}$. Thus the closed set $\overline{S}$ is also open in $\hat{X}$, and so also the first map is well-defined.

If $S \in \mathcal{OC}_X$, then $S \subseteq \overline{S} \cap X$. Let $x \in \overline{S} \cap X$ and suppose that $x \in X \backslash S$. Then $x \in \overline{X \backslash S} = \hat{X} \backslash \overline{S}$ which is a contradiction, and it follows that $x \in S$. Therefore $S = \overline{S} \cap X$.

If $\mathcal{S} \in \mathcal{OC}_{\hat{X}}$, then $\overline{\mathcal{S} \cap X} \subseteq \mathcal{S}$, since $\mathcal{S}$ is closed. Since $\mathcal{S}$ is also open, $\mathcal{S} \cap X$ is dense in $\mathcal{S}$, and so $\overline{\mathcal{S} \cap X} = \mathcal{S}$. This proves that the considered maps are bijections.

In order to prove (ii) let $\hat{x} \in \hat{X}_{iso}$. Then $\{\hat{x}\}$ is open in $\hat{X}$. Since $i(X)$ is dense in $\hat{X}$, the set $\{\hat{x}\} \cap i(X)$ is not empty and so $\hat{x} \in i(X)$. Thus $\{\hat{x}\}$ is an open subset of $i(X)$.

Conversely, let $x \in X_{iso}$. Since the set $\{x\}$ is open and closed in $X$, the same is true, by (i), for its closure $\overline{\{x\}}$ in $\hat{X}$. Consider the open set

$$U = \overline{\{x\}} \backslash \{i(x)\} \subseteq \hat{X}$$

(observe that $\{i(x)\}$ is closed in the Hausdorff space $\hat{X}$). Suppose that $U$ is not empty. Then, using (i), we get the contradiction

$$\emptyset \neq U \cap i(X) = (\overline{\{x\}} \cap i(X)) \backslash \{i(x)\} = \{i(x)\} \backslash \{i(x)\}.$$

Therefore $U$ is empty, i.e. $\overline{\{x\}} = \{i(x)\}$, and so $\{i(x)\}$ is open in $\hat{X}$. □

Now let $(X, \mathcal{T}) = (\mathcal{P}_K, \mathcal{T}_K)$. This space is a strongly zero-dimensional Hausdorff space by proposition 3.5(ii). If $\mathfrak{U}_K = \mathfrak{U}_K^{oc}$ denotes the uniformity of finite partitions of $\mathcal{P}_K$ by both open and closed subsets of $(\mathcal{P}_K, \mathcal{T}_K)$, then we obtain

**Theorem 4.4** *The Hausdorff uniform space $(\mathcal{P}_K, \mathfrak{U}_K)$ is pre-compact and strongly zero-dimensional, and its completion $(\hat{\mathcal{P}}_K, \hat{\mathfrak{U}}_K)$ is a profinite space. The canonical map*

$$i : (\mathcal{P}_K, \mathfrak{U}_K) \longrightarrow (\hat{\mathcal{P}}_K, \hat{\mathfrak{U}}_K)$$

*induces an isomorphism of $(\mathcal{P}_K, \mathfrak{U}_K)$ onto a dense subspace of $(\hat{\mathcal{P}}_K, \hat{\mathfrak{U}}_K)$.*

*Furthermore the sets $(\mathcal{P}_K)_{iso}$ and $(\hat{\mathcal{P}}_K)_{iso}$ of isolated points are isomorphic and finite.*

4. *The uniformity of finite partitions by both open and closed and finitely presented subsets of $(\mathcal{P}_K, \mathcal{T}_K)$.*

We define

$$Part^{ocf} = \{(V_1, \ldots, V_n) \mid n \in \mathbb{N}, V_i \text{ open, closed, finitely presented}, \dot{\bigcup}_{i=1}^n V_i = \mathcal{P}_K\},$$

i.e. this is the set of finite partitions of $\mathcal{P}_K$ by open and closed subsets of $\mathcal{T}_K$ which are finitely presented. Let

$$\mathfrak{V}_K^{ocf} = \{V_Q = \dot{\bigcup}_{i=1}^n (V_i \times V_i) \subseteq \mathcal{P}_K \times \mathcal{P}_K \mid n \in \mathbb{N}, Q = (V_1, \ldots, V_n) \in Part^{ocf}\}.$$

Then $\mathfrak{V}_K^{ocf}$ is a base for a uniform structure on $\mathcal{P}_K$ denoted by $\mathfrak{U}_K^{ocf}$. Obviously, $\mathfrak{U}_K^{ocf}$ is coarser than $\mathfrak{U}_K = \mathfrak{U}_K^{oc}$ (and it seems unlikely that they are equal). By lemma 3.7 the uniformity $\mathfrak{U}_K^{ocf}$ induces the topology $\mathcal{T}_K$ on $\mathcal{P}_K$.

5. *The uniformity of finite partitions by open and closed subsets of $(\mathcal{P}_K, \mathcal{T}_K)$ defined by the discriminant of finite Galois extensions $F|E$, $E \subseteq K$.*

Let $d \in \mathbb{N}$ and let

$$S_d = S_d(K) = \{F|E \text{ a finite Galois extension}, E \subseteq K, |D(F|\mathbb{Q})| = d\},$$

where $D(F|\mathbb{Q})$ denotes the discriminant of $F$. The set $S_d$ is finite by Hermite's theorem (and can be empty), see [3] III. (2.16). For $x \in \mathcal{P}_K$ let

$$S_{d,x} = S_{d,x}(K) = \{F|E \text{ finite Galois }, E \subseteq K, x \in U^K(F|E), |D(F|\mathbb{Q})| = d\}$$

and

$$V_d(x) = \begin{cases} \bigcap_{F|E \in S_{d,x}} P_{F|E}^K\left(\left(\dfrac{F|E}{x_F}\right)\right), & \text{if } S_{d,x} \neq \varnothing, \\ \mathcal{P}_K, & \text{otherwise,} \end{cases}$$

where $x_F$ is an extension of $x_E = E \cap x$ to $F$. Furthermore let

$$R_d = R_d(K) = \bigcup_{F|E \in S_d} R^K(F|E),$$

$$G_d = G_d(K) = \prod_{F|E \in S_d} G(F|E).$$

The elements $(\sigma_{F|E})_{F|E}$ of $G_d$ will be denoted by $\tilde{\sigma}$.

For an open finitely presented set $U$ let

$$d(U) = \inf\{\max_{ij}|D(F_{ij}|\mathbb{Q})|, \text{ if } U = \bigcup_{i=1}^{n}\bigcap_{j=1}^{n_i} P^K_{F_{ij}|E_{ij}}(\sigma_{ij})\},$$

where the infimum is taken over all possible finite presentations of $U$. Assume that $S_d \neq \varnothing$. By lemma 3.7 we find pairwise disjoint, open and closed neighbourhoods $U(\alpha)$ of the elements $\alpha \in R_d$ which are finitely presented and contained in $\bigcap_{m \leq d} V_m(\alpha)$. Let

$$d(R_d) = \inf\{\max\{d(U(\alpha)), \alpha \in R_d\}\},$$

where the infimum is taken over the set of all these coverings of $R_d$. The set $C^{\min}$ of coverings $(U(\alpha)), \alpha \in R_d)$ with $\max\{d(U(\alpha)), \alpha \in R_d\} = d(R_d)$ is finite since there are only finitely many fields $F$ with $|D(F|\mathbb{Q})| \leq d(R_d)$. It follows that the covering

$$(\tilde{V}_d(\alpha), \alpha \in R_d) = \bigcap_{C^{\min}}(U(\alpha), \alpha \in R_d),$$

of $R_d$ is given by pairwise disjoint, open and closed neighbourhoods $\tilde{V}_d(\alpha)$, which are contained in $\bigcap_{m \leq d} V_m(\alpha)$, and this covering is uniquely determined by $d$.

Now we define $V_d = \mathcal{P}_K \times \mathcal{P}_K$ if $S_d = \varnothing$, and otherwise

$$V_d = \Big((V_d(0)\backslash V_{R_d}) \times (V_d(0)\backslash V_{R_d})\Big) \cup \bigcup_{\alpha \in R_d}\Big(\tilde{V}_d(\alpha) \times \tilde{V}_d(\alpha)\Big),$$

where

$$V_d(0) = \bigcup_{\tilde{\sigma} \in G_d}\bigcap_{F|E \in S_d}\Big(P^K_{F|E}(\sigma_{F|E})\Big), \qquad V_{R_d} = \bigcup_{\alpha \in R_d}\tilde{V}_d(\alpha)$$

and

$$V_d(0)\backslash V_{R_d} = \bigcup_{\tilde{\sigma} \in G_d}\bigcap_{F|E \in S_d}\Big(P^K_{F|E}(\sigma_{F|E})\backslash V_{R_d}\Big).$$

If $S_d \neq \varnothing$ and $x \in \mathcal{P}_K$ is not contained in $R_d$, then $x \in V_d(0)$. Therefore

$$(V_d(0)\backslash V_{R_d}) \cup \bigcup_{\alpha \in R_d}\tilde{V}_d(\alpha)$$

is a partition of $\mathcal{P}_K$ by open and closed sets which are finitely presented. We put

$$W_d = \bigcap_{m \leq d} V_m.$$

Then $W_{d'} \subseteq W_d$ for $d \leq d'$ and $\mathfrak{V}^D_K = \{W_d, d \in \mathbb{N}\}$ is a base for a uniform structure on $\mathcal{P}_K$ denoted by $\mathfrak{U}^D_K$.

**Lemma 4.5** *The uniformity $\mathfrak{U}_K^D$ is equal to $\mathfrak{U}_K^{ocf}$.*

**Proof:** Obviously, the uniformity $\mathfrak{U}_K^D$ is coarser than $\mathfrak{U}_K^{ocf}$. Let

$$\bigcup_{i=1}^{n} (V_i \times V_i) \in \mathfrak{V}_K^{ocf} \quad \text{and} \quad d = \max\{d(V_i),\, i=1,\ldots,n\}.$$

We claim that

$$\bigcap_{m \leq d} \left( (V_m(0) \setminus V_{R_m}) \cup \bigcup_{\alpha \in R_m} \tilde{V}_m(\alpha) \right) = \mathcal{P}_K \quad \text{is a refinement of} \quad \bigcup_{i=1}^{n} V_i = \mathcal{P}_K.$$

Indeed, every set of the first partition is contained in an open set of the form $\bigcap_{m \leq d} V_m(x_0)$ for some $x_0 \in \mathcal{P}_K$. Let $i_0$ be the number such that $x_0 \in V_{i_0}$ and let

$$V_{i_0} = \bigcup_{k=1}^{r} \bigcap_{j=1}^{r_k} P_{F_{kj}|E_{kj}}^K(\sigma_{kj})$$

be a presentation of $V_{i_0}$ such that $d(V_{i_0}) = \max_{kj} |D(F_{kj}|\mathbb{Q})|$. It follows that $x_0 \in \bigcap_{j=1}^{r_k} P_{F_{kj}|E_{kj}}^K(\sigma_{kj})$ for some $k$. Since $\max_j |D(F_{kj}|\mathbb{Q})| \leq d(V_{i_0}) \leq d$ we get $\bigcap_{m \leq d} V_m(x_0) \subseteq V_{i_0}$. This proves the claim, and so $\mathfrak{U}_K^D$ is finer than $\mathfrak{U}_K^{ocf}$. $\square$

# 5 Metric

In this section we will define a metric on $\mathcal{P}_K$ which induces the uniformity $\mathfrak{U}_K^D$ and therefore the topology $\mathcal{T}_K$. The idea is that two points $x, y \in \mathcal{P}_K$ are *near*, if they induce in *many* fields with *large* discriminant the same Frobenius automorphism.

**Theorem 5.1** *The map*

$$\delta : \mathcal{P}_K \times \mathcal{P}_K \longrightarrow [0,1], \quad (x,y) \mapsto \delta(x,y) = \frac{1}{n},$$

*where*

$$n = \sup\{\, d \mid (x,y) \in W_d\},$$

*defines an ultra-metric on $\mathcal{P}_K$ which induces the uniformity $\mathfrak{U}_K^D$.*

**Proof:** Obviously, $\delta$ is symmetric, $\delta(x,x) = 0$ for every $x \in \mathcal{P}_K$,

$$\delta(x,y) \leq \max(\delta(x,z), \delta(z,y)) \quad \text{for all } x, y, z \in \mathcal{P}_K$$

and the (quasi-)metric $\delta$ induces the uniformity $\mathfrak{U}_K^D$. Since $\mathcal{T}_K$ is a Hausdorff topology, $\delta$ is an (ultra-)metric. $\square$

**Corollary 5.2** *The completion $(\hat{\mathcal{P}}_K, \hat{\mathfrak{U}}_K^D)$ of $(\mathcal{P}_K, \mathfrak{U}_K^D)$ is a profinite space.*

**Proof:** Since $\mathfrak{U}_K^D$ is coarser than $\mathfrak{U}_K$, we get a surjection

$$(\hat{\mathcal{P}}_K, \hat{\mathfrak{U}}_K) \longrightarrow (\hat{\mathcal{P}}_K, \hat{\mathfrak{U}}_K^D),$$

and so $(\hat{\mathcal{P}}_K, \hat{\mathfrak{U}}_K^D)$ is compact. Furthermore, the extension of $\delta$ to $(\hat{\mathcal{P}}_K, \hat{\mathfrak{U}}_K^D)$ is also an ultra-metric and so the completion is totally disconnected. $\square$

It is obvious that the metric $\delta$ is not easily to calculate even in the case $K = \mathbb{Q}$ (at least if $d > 21$ when not only quadratic fields are involved). But we do some calculations for $d \leq 4$. We have two extensions $F|\mathbb{Q}$ with absolute discriminant $|D(F|\mathbb{Q})| \leq 4$:

$$\mathbb{Q}(\sqrt{-3}),\ |D| = 3, \quad \text{and} \quad \mathbb{Q}(\sqrt{-1}),\ |D| = 4,$$

but we will also need the field $\mathbb{Q}(\sqrt{5})$ with $|D| = 5$. From now on we use the notation $(a|\gamma)$ for $P_{\mathbb{Q}(\sqrt{a})|\mathbb{Q}}(\gamma)$, $\gamma \in G(\mathbb{Q}(\sqrt{a})|\mathbb{Q})$, and we denote the non-trivial element of $G(\mathbb{Q}(\sqrt{a})|\mathbb{Q})$ by $-1$. Then

$$V_3(0) = (-3|\,1) \cup (-3|-1), \quad \tilde{V}_3(3) = (-1|-1) \cap (5\,|-1),$$

and the partition of $\mathcal{P}_\mathbb{Q}$ by both open and closed sets associated to $W_3$ is

$$\begin{aligned}\mathcal{P}_\mathbb{Q} &= (-3|\,1)\backslash\Big((-1|-1) \cap (5\,|-1)\Big) \cup (-3|-1)\backslash\Big((-1|-1) \cap (5\,|-1)\Big) \\ &\quad \cup \Big((-1|-1) \cap (5\,|-1)\Big).\end{aligned}$$

For $d = 4$ we get

$$V_4(0) = (-1|\,1) \cup (-1|-1), \quad \tilde{V}_4(2) = (-3|\,1) \cap (5\,|\,1),$$

and, calculating the corresponding partition associated to $W_4$, we obtain

$$\begin{aligned}
\mathcal{P}_\mathbb{Q} &= \Big((-3|\,1)\cap(-1|\,1)\Big)\backslash\Big(\tilde{V}_3(3)\cup\tilde{V}_4(2)\Big) \\
&\cup \Big((-3|\,1)\cap(-1|-1)\Big)\backslash\Big(\tilde{V}_3(3)\cup\tilde{V}_4(2)\Big) \\
&\cup \Big((-3|-1)\cap(-1|\,1)\Big)\backslash\Big(\tilde{V}_3(3)\cup\tilde{V}_4(2)\Big) \\
&\cup \Big((-3|-1)\cap(-1|-1)\Big)\backslash\Big(\tilde{V}_3(3)\cup\tilde{V}_4(2)\Big) \\
&\cup \tilde{V}_3(3)\backslash\tilde{V}_4(2) \\
&\cup \tilde{V}_4(2)\backslash\tilde{V}_3(3).
\end{aligned}$$

It follows that for different prime numbers $x, y \leq 19$

$\delta(x, y) = \frac{1}{3}$, if $(x, y), (y, x) = (2, 7), (2, 13), (5, 11), (7, 13), (7, 19), (13, 19)$,

$\delta(x, y) \leq \frac{1}{4}$, if $(x, y), (y, x) = (2, 19), (5, 17)$,

and for all other pairs we have $\delta(x, y) = 1$.

**Concluding remark:** For $m \in \mathbb{N}$, $m \geq 2$, we define

$$\mathcal{B}_K^m = \{P_{F|E}^K(\sigma) \,|\, E \subseteq K, F|E \text{ a finite Galois extension of degree } [F:E] \leq m,$$
$$\sigma \in G(F|E)\}.$$

Let $\mathcal{T}_K^m$ be the topology on $\mathcal{P}_K$ with subbase $\mathcal{B}_K^m$. Then every statement of this and the preceding section remains true. Crucial is proposition 3.4(i) and its proof is not infected if $m > 2$. For $m = 2$ we need a different argument.

**Proposition 5.3** *For every two different points $\mathfrak{p}_1$ and $\mathfrak{p}_2$ of $(\mathcal{P}_K, \mathcal{T}_K^2)$ there exist both open and closed neighbourhoods $W_1 \in \mathcal{T}_K^2$ and $W_2 \in \mathcal{T}_K^2$ of $\mathfrak{p}_1$ and $\mathfrak{p}_2$, respectively, such that*
$$W_1 \cap W_2 = \emptyset.$$

**Proof:** First we reduce to the case $K = \mathbb{Q}$. By the theorem of Grunwald/Wang there exist *quadratic* extensions $L_i|K$, $i = 1, 2$, such that $\mathfrak{p}_1$ and all its conjugates are inert and $\mathfrak{p}_2$ is completely decomposed in $L_1|K$ and the same holds if we interchange the indices 1 and 2. If $\sigma_i$ denotes the non-trivial element of $G(L_i|K)$, then it follows that

$$\mathfrak{p}_1 \in P_{L_1|K}(\sigma_1), \mathfrak{p}_2 \notin P_{L_1|K}(\sigma_1) \text{ and } \mathfrak{p}_2 \in P_{L_2|K}(\sigma_2), \mathfrak{p}_1 \notin P_{L_2|K}(\sigma_2),$$

and
$$R(L_i|K) \cap \{\mathfrak{p}_i' \mid \mathfrak{p}_i' \text{ conjugated to } \mathfrak{p}_i\} = \varnothing, \quad i = 1, 2.$$

Let $U_i = P_{L_i|K}(\sigma_i)$, $i = 1, 2$. These open neighbourhoods need not to be disjoint or closed. Let
$$\varphi_{K|\mathbb{Q}}(R(L_1|K)) = \{q_1, \ldots, q_s\}$$
and let $p_1$ be the underlying prime number of $\mathfrak{p}_1$ which is by construction different to the prime numbers $q_1, \ldots, q_s$. Assume that we have proved the proposition in the case $K = \mathbb{Q}$. Then, using the argument of the proof of proposition 3.4(ii), there are pairwise disjoint both open and closed neighbourhoods $U(p_1), U(q_1), \ldots, U(q_s) \in \mathcal{T}_\mathbb{Q}^2$ of the considered prime numbers. It follows that

$$\mathfrak{p}_1 \in U_1 \cap \varphi_{K|\mathbb{Q}}^{-1}(U(p_1)) \quad \text{and} \quad R(L_1|K) \subseteq \bigcup_{i=1}^s \varphi_{K|\mathbb{Q}}^{-1}(U(q_i)).$$

The open set $V_1 = U_1 \cap \varphi_{K|\mathbb{Q}}^{-1}(U(p_1))$ is also closed since

$$\begin{aligned}
\overline{U_1 \cap \varphi_{K|\mathbb{Q}}^{-1}(U(p_1))} &\subseteq \overline{U_1} \cap \varphi_{K|\mathbb{Q}}^{-1}(U(p_1)) \\
&\subseteq (U_1 \cup R(L_1|K)) \cap \varphi_{K|\mathbb{Q}}^{-1}(U(p_1)) \\
&= U_1 \cap \varphi_{K|\mathbb{Q}}^{-1}(U(p_1))
\end{aligned}$$

(observe that the map $\varphi_{K|\mathbb{Q}}$ is continuous and that the sets $\varphi_{K|\mathbb{Q}}^{-1}(U(p_1))$ and $\varphi_{K|\mathbb{Q}}^{-1}(U(q_i))$, $i = 1, \ldots, s$, are also pairwise disjoint). Constructing an open and closed neighbourhood $V_2$ of $\mathfrak{p}_2$ in the same way and observing that

$$\mathfrak{p}_1 \in V_1, \mathfrak{p}_2 \notin V_1 \text{ and } \mathfrak{p}_2 \in V_2, \mathfrak{p}_1 \notin V_2,$$

then the neighbourhoods $W_1 = V_1 \backslash V_2$ of $\mathfrak{p}_1$ and $W_2 = V_2 \backslash V_1$ of $\mathfrak{p}_2$ have the desired properties.

Now let $K = \mathbb{Q}$ and let $p_1$ and $p_2$ be two different prime numbers. Assume that $p_1$ is odd and let
$$\sigma = \left(\frac{p_2}{p_1}\right) \in \{\pm 1\}.$$

Let $q$ be an odd prime number, different to $p_1$ and $p_2$, such that
$$\left(\frac{p_1}{q}\right) = -\sigma, \quad \left(\frac{p_2}{q}\right) = \sigma, \quad \left(\frac{-1}{q}\right) = 1.$$

Such a prime number exists since the set
$$P_{\mathbb{Q}(\sqrt{p_1})|\mathbb{Q}}(-\sigma) \cap P_{\mathbb{Q}(\sqrt{p_2})|\mathbb{Q}}(\sigma) \cap P_{\mathbb{Q}(\sqrt{-1})|\mathbb{Q}}(1),$$
has positive density, and so it is not empty. Let
$$L_q = \mathbb{Q}(\sqrt{(-1)^{\varepsilon(q)}q}) \quad \text{and} \quad L_{p_1} = \mathbb{Q}(\sqrt{(-1)^{\varepsilon(p_1)}p_1}),$$

where $\varepsilon(\ell) = (\ell-1)/2$. Then $p_1 \notin P_{L_q|\mathbb{Q}}(\sigma)$ and $q \notin P_{L_{p_1}|\mathbb{Q}}(\sigma)$ since

$$\left(\frac{(-1)^{\varepsilon(q)}q}{p_1}\right) = \left(\frac{p_1}{q}\right) = -\sigma \quad \text{and} \quad \left(\frac{(-1)^{\varepsilon(p_1)}p_1}{q}\right) = \left(\frac{p_1}{q}\right)\left(\frac{-1}{q}\right)^{\varepsilon(p_1)} = -\sigma.$$

Furthermore, $p_2 \in P_{L_q|\mathbb{Q}}(\sigma) \cap P_{L_{p_1}|\mathbb{Q}}(\sigma)$: if $p_2$ is odd, then

$$\left(\frac{(-1)^{\varepsilon(q)}q}{p_2}\right) = \left(\frac{p_2}{q}\right) = \sigma \quad \text{and} \quad \left(\frac{(-1)^{\varepsilon(p_1)}p_1}{p_2}\right) = \left(\frac{p_2}{p_1}\right) = \sigma.$$

If $p_2 = 2$, then $p_2$ splits in $L_q|\mathbb{Q}$ if and only if $(-1)^{\varepsilon(q)}q \equiv 1 \mod 8$, i.e. if and only if $\left(\frac{2}{q}\right) = 1$. Since $\sigma = \left(\frac{2}{q}\right) = \left(\frac{2}{p_1}\right)$, we have

$$2 \in P_{L_q|\mathbb{Q}}(\sigma),$$

and also $2 \in P_{L_{p_1}|\mathbb{Q}}(\sigma)$.

It follows that $V_2 = P_{L_{p_1}|\mathbb{Q}}(\sigma) \cap P_{L_q|\mathbb{Q}}(\sigma)$ is a neighbourhood of $p_2$ which is both open and closed (since $R(L_q|\mathbb{Q}) = \{q\}$ and $R(L_{p_1}|\mathbb{Q}) = \{p_1\}$), and does not contain $p_1$.

If $p_2$ is odd, then we construct an open and closed neighbourhood $V_1$ of $p_1$ in the same way, and the neighbourhoods $W_1 = V_1 \backslash V_2$ of $p_1$ and $W_2 = V_2 \backslash V_1$ of $p_2$ have the desired properties.

If $p_2 = 2$, then we have also to construct an open and closed neighbourhood $V_1$ of $p_1$ with the analogous properties. Again let

$$\sigma = \left(\frac{2}{p_1}\right),$$

and let $q$ be an odd prime number such that

$$\left(\frac{p_1}{q}\right) = \sigma \quad \text{and} \quad \left(\frac{2}{q}\right) = -\sigma.$$

If

$$L_q = \mathbb{Q}(\sqrt{(-1)^{\varepsilon(q)}q}) \quad \text{and} \quad L_2 = \mathbb{Q}(\sqrt{2}),$$

then $p_1 \in P_{L_q|\mathbb{Q}}(\sigma) \cap P_{L_2|\mathbb{Q}}(\sigma)$:

$$\left(\frac{(-1)^{\varepsilon(q)}q}{p_1}\right) = \left(\frac{p_1}{q}\right) = \sigma, \quad \left(\frac{2}{p_1}\right) = \sigma.$$

Furthermore, $q \notin P_{L_2|\mathbb{Q}}(\sigma)$ because

$$\left(\frac{2}{q}\right) = -\sigma,$$

and $2 \notin P_{L_q|\mathbb{Q}}(\sigma)$ because $2 \in P_{L_q|\mathbb{Q}}(\left(\frac{2}{q}\right)) = P_{L_q|\mathbb{Q}}(-\sigma)$. It follows that $V_1 = P_{L_2|\mathbb{Q}}(\sigma) \cap P_{L_q|\mathbb{Q}}(\sigma)$ is a neighbourhood of $p_1$ which is both open and closed and does not contain $p_2$. This finishes the proof of the proposition.

Mathematisches Institut
der Universität Heidelberg
Im Neuenheimer Feld 288
69120 Heidelberg
Germany

e-mail: wingberg@mathi.uni-heidelberg.de